\documentclass[12pt]{article}
\usepackage[utf8]{inputenc}
\usepackage[english]{babel}
\usepackage[pdftex]{color}
\usepackage[pdftex]{graphicx}
\usepackage{amsmath,amssymb,amsthm}
\usepackage{tikz}
\usepackage{caption}
\usepackage{subcaption}
\usepackage{booktabs}

\usepackage[a4paper,left=15mm,right=15mm,top=30mm,bottom=20mm]{geometry}
\newtheorem{theor_eng}{Theorem}[section]
\newtheorem{lemma}[theor_eng]{Lemma}
\newtheorem{corollary}[theor_eng]{Corollary}
\newtheorem{remark}[theor_eng]{Remark}

\begin{document}

\bibliographystyle{unsrt}

\title{Balanced affine Motzkin paths: Pearson geometry and global endpoint asymptotics}

\author{
Alexander Omelchenko\\
\small Constructor University Bremen, Campus Ring 1, 28759 Bremen, Germany\\
\small\tt aomelchenko@constructor.university
}

\maketitle

\begin{abstract}
We study endpoint distributions of balanced affine weighted Motzkin
paths. In the balanced case, the generating-function equation has
Pearson-type characteristic geometry. We show that this geometry controls
the terminal-height law globally: the characteristic escape time
determines the limiting cumulant generating function, the
large-deviation rate function, and the ray-scale asymptotics. Thus the
usual Gaussian window is only the local quadratic approximation to a
global Pearson-driven profile. For finite sizes, we prove a uniform
Daniels saddlepoint approximation in the one-dominant-singularity
regimes and identify the exceptional antipodal case requiring a
lattice/interference correction.

\end{abstract}

\section*{Introduction}

Lattice paths are classical and central objects of enumerative
combinatorics \cite{Stanley2,Flajolet2009,Krattenthaler2015}. A basic
example is the \emph{Motzkin path}: a walk on the integer lattice that
starts at the origin, uses up-, level-, and down-steps, and never
descends below the horizontal axis. The corresponding \emph{Motzkin
triangle} records such paths by length and terminal height.

Weighted versions of this triangle occupy a particularly important
place in the subject. Through the Flajolet--Viennot correspondence,
height-weighted Motzkin paths encode continued fractions, three-term
recurrences, and the moment theory of orthogonal polynomials
\cite{Flajolet1980CF,Viennot1985,Chihara1978,Ismail2005}. In this
interpretation, choosing the step weights amounts to choosing Jacobi
data: the bottom row of the triangle gives the corresponding moment
sequence, while the full triangle refines it by terminal height.

\begin{figure}[ht]
\centering
\includegraphics{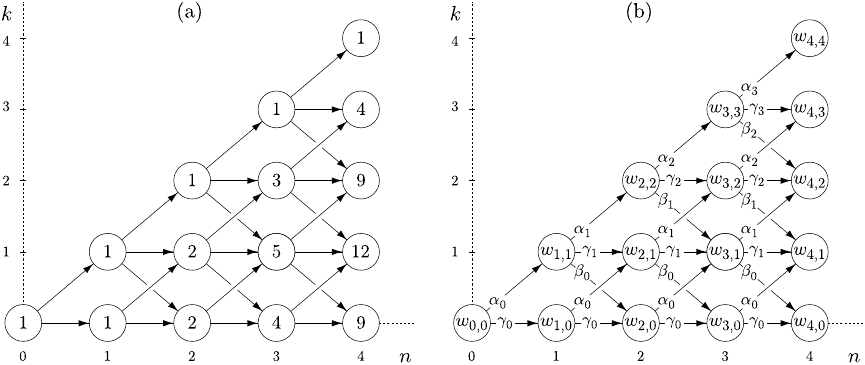}
\caption{The classical Motzkin triangle (a) and a triangle with linearly
varying multiplicities (b).}
\label{fig2}
\end{figure}

The present paper studies the affine height-dependent version of this
model, in which the step weights are
\[
 \alpha_k = a\,k+\alpha_0,\qquad
 \beta_k  = b\,k+\beta_0,\qquad
 \gamma_k = c\,k+\gamma_0 .
\]
Here $\alpha_k$ and $\gamma_k$ weight the up- and level-steps leaving
height $k$, while $\beta_k$ weights the down-step arriving at height
$k$. Throughout we write $A=a$, $B=c$, $C=b$,
$Q(x)=Ax^{2}+Bx+C$, and $\Delta=B^{2}-4AC$. For each $n$ and $k$, the
number $w_{n,k}$ denotes the total weight of paths of length $n$ that
terminate at height $k$, and $P_n(x)=\sum_{k}w_{n,k}x^{k}$ is the row
polynomial of the triangle.

The affine Motzkin and Dyck triangles were introduced and studied
enumeratively in \cite{Meshkov2010}. The aim of the present paper is
different. We do not seek another enumeration of the affine triangles.
Instead we study the law of the terminal height of long paths: a path
of length $n$ is drawn with probability proportional to its weight, and
the resulting distribution is $p_{n,k}=w_{n,k}/P_n(1)$, the law of the
random terminal height $K_n$.
From the orthogonal-polynomial viewpoint, the bottom row of the
triangle is only a moment sequence; what we analyse is the full
terminal-height refinement of that moment sequence.
Our main point is that, in the
\emph{balanced} affine model, the entire asymptotic profile of this law
is controlled by a Pearson-type characteristic geometry, through the
chain
\[
\begin{array}{c}
\text{balanced affine Motzkin recurrence}\\
\Downarrow\\
\text{local Pearson-type PDE}\\
\Downarrow\\
\text{moving singularity } t=\tau(x)\\
\Downarrow\\
F(\theta)=\log\dfrac{\tau(1)}{\tau(e^{\theta})}\\
\Downarrow\\
I(u)=\sup_{\theta\in\mathbb R}\{u\theta-F(\theta)\}\\
\Downarrow\\
p_{n,k}\ \text{ on the ray scale } k\sim un .
\end{array}
\]
Here $\tau(x)=\int_x^\infty dy/Q(y)$ is the escape time of the
characteristic flow. In the quadratic regimes $A>0$ it is finite, and it
coincides with the first positive singularity, in the $t$-plane, of the
exponential generating function
$w(x,t)=\sum_{n\ge0}P_n(x)\,t^{n}/n!$.

The connection with Pearson distributions is structural rather than
literal: the terminal-height law is not itself a Pearson distribution.
The Pearson feature is the logarithmic derivative along
characteristics. In the balanced case the generating function satisfies
a local first-order PDE, and eliminating the characteristic time gives
\[
  \frac{d}{dx}\log w
  =
  -\,\frac{\alpha_0x+\gamma_0}{Ax^{2}+Bx+C},
\]
a linear numerator over a quadratic denominator. This is exactly the
form that organises the classical Pearson system of probability
distributions \cite{Pearson1895,JohnsonKotz1994}. The same quadratic
$Q$ drives the characteristic flow, determines the escape time
$\tau(x)$, and ultimately governs the rate function $I(u)$.

This was the motivating observation behind the paper: a Pearson-type
equation appears, unexpectedly, in a weighted lattice-path problem. The
question is whether this occurrence is a formal analogy, or whether it
actually controls the probability distribution associated with the
paths. The answer developed below is that it controls the full endpoint
profile---locally through a Gaussian window, globally through the
large-deviation rate function, and at finite $n$ through a saddlepoint
approximation.

The balanced condition is the structural reason why the mechanism is
visible. For the general six-parameter affine family, the
generating-function PDE contains a nonlocal boundary term proportional
to $\beta_0-C$ and evaluated at $x=0$. The condition $\beta_0=C$ is
precisely the condition under which this term disappears. The balanced
models therefore form the local, solvable core of the affine family:
the balanced condition is the locality condition of the PDE, not an ad
hoc restriction.

The affine class itself sits at a natural interface, and this is the
broader reason for studying it. Under the Flajolet--Viennot
correspondence, affine step weights produce Jacobi data with affine
diagonal part and quadratic off-diagonal products---the pattern
underlying the classical Meixner--Sheffer orthogonal families
\cite{Meixner1934,Chihara1978,Ismail2005}. Closely related affine
tridiagonal operators occur as generators of birth--death
processes \cite{KarlinMcGregor1957,Anderson1991}---a basic language for
stochastic population and reaction models---and as scalar analogues of
level-dependent quasi-birth--death queues
\cite{Neuts1981,LatoucheRamaswami1999}. In random matrix theory,
tridiagonal and Jacobi models provide sparse realisations of beta
ensembles \cite{DumitriuEdelman2002,KillipNenciu2004,Forrester2010}. In
statistical physics, Motzkin paths appear verbatim as the ground-state
basis of exactly solvable quantum spin chains, where the height
distribution of a random weighted Motzkin path at a cut governs
entanglement and correlation properties
\cite{BravyiCMNS2012,MovassaghShor2016}; the weights arising there
(colour multiplicities, area deformations) differ from the affine
height weights considered here, but the question is of the same type:
asymptotic statistics of a height functional of a weighted Motzkin
ensemble.

In all these settings, what is typically needed beyond moment or
stationary information is precisely endpoint and rare-event behaviour:
tail probabilities, finite-horizon profiles, finite-size corrections.
The present paper does not claim a direct application to any of these
areas. Rather, it isolates an exactly solvable scalar model in which
the tridiagonal Pearson structure is explicit enough to yield not only
closed forms and moments, but the full endpoint probability profile:
large deviations, ray-scale prefactors, and uniform finite-$n$
saddlepoint approximations.

The Pearson discriminant is also visible combinatorially. Several
classical moment sequences occur as bottom rows $w_{n,0}$ of balanced
affine triangles; Table~\ref{tab:classical} lists one representative
from each of the five regimes. In particular, among the genuinely
quadratic rows, the same discriminant that classifies the
characteristic equation separates permutations ($\Delta=0$), even
alternating permutations ($\Delta<0$) \cite{Andre1879}, and ordered set
partitions ($\Delta>0$). The discriminant is therefore not only an
analytic parameter; it is reflected in classical combinatorial
families.

\begin{table}[ht]
\centering
\scriptsize
\resizebox{\textwidth}{!}{%
\begin{tabular}{@{}llllll@{}}
\toprule
regime & $(\alpha_k,\beta_k,\gamma_k)$ & $w_{n,0}$ & $w(0,t)$ & objects & OPS family\\
\midrule
$A=B=0$
  & $(1,k{+}1,0)$
  & $1,0,1,0,3,\dots$
  & $e^{t^2/2}$
  & perfect matchings
  & Hermite\\
$A=0,\ B>0$
  & $(1,k{+}1,k{+}1)$
  & $1,1,2,5,15,\dots$
  & $e^{e^t-1}$
  & set partitions
  & Charlier\\
$A>0,\ \Delta=0$
  & $(k{+}1,k{+}1,2k{+}1)$
  & $n!$
  & $(1-t)^{-1}$
  & permutations
  & Laguerre $(\alpha=0)$\\
$A>0,\ \Delta<0$
  & $(k{+}1,k{+}1,0)$
  & $1,0,1,0,5,\dots$
  & $\sec t$
  & even alternating permutations
  & Meixner--Pollaczek\\
$A>0,\ \Delta>0$
  & $(k{+}1,2(k{+}1),3k{+}1)$
  & $1,1,3,13,75,\dots$
  & $(2-e^t)^{-1}$
  & ordered set partitions
  & Meixner $(\beta=1,c=1/2)$\\
\bottomrule
\end{tabular}%
}
\caption{Classical specialisations inside the balanced affine Motzkin
model. In each case the bottom row $w_{n,0}$ is a classical moment
sequence, and the full triangle refines it by terminal height. In the
usual Motzkin-path normalisation the corresponding Jacobi data are
$b_h=\gamma_h$ and $\lambda_h=\alpha_{h-1}\beta_{h-1}$. All five
specialisations satisfy the balanced condition $\beta_0=C$, introduced
below. The OPS identifications are meant up to the standard
normalisation of monic three-term recurrences.}
\label{tab:classical}
\end{table}

We now describe the asymptotic results. A first analysis of the
coefficients yields the expected local picture: the law of $K_n$
concentrates around a mean of order $n$, with Gaussian fluctuations in
a $\sqrt n$-window in the one-dominant-singularity regimes, and with
analogous statements at sublinear scales in the degenerate regimes
$A=0$. This local statement, however, is too coarse to reveal the
Pearson structure: near the mean, all regimes look Gaussian to first
order. The distinctions between the Pearson regimes appear only on the
ray scale $k\sim un$, where the whole movement of the singularity
$\tau(x)$ under exponential tilting becomes relevant.

We therefore pass to exponential tilting. Writing $x=e^{\theta}$ and
$\kappa_n(\theta)=\log P_n(e^{\theta})$, the tilt reweights the law so
that any prescribed height becomes typical: for a target $k$, the
saddlepoint $\theta_{n,k}$ is the unique tilt at which $k$ equals the
mean of the tilted distribution, $\kappa_n'(\theta_{n,k})=k$. Expanding
the lattice Fourier inversion integral at this saddle gives Daniels'
saddlepoint approximation \cite{Daniels1954}: a single finite-$n$
formula that agrees with the Gaussian approximation at the centre and
remains pointwise accurate deep into the tails.

On the $n$-scale, in the quadratic regimes, the same tilted geometry
produces the limit cumulant generating function
\[
  F(\theta)
  =
  \lim_{n\to\infty}\frac1n\log\frac{P_n(e^{\theta})}{P_n(1)}
  =
  \log\frac{\tau(1)}{\tau(e^{\theta})}.
\]
In the non-degenerate case ($A>0$, $\alpha_0>0$, $B+C>0$) this function
is smooth and strictly convex, and its derivative maps $\mathbb R$
bijectively onto $(0,1)$; the finite-$n$ saddlepoint equation
$\kappa_n'(\theta)=k$ is the finite-$n$ counterpart of the limiting
equation $F'(\theta)=u$. The Legendre transform
\[
  I(u)=\sup_{\theta\in\mathbb R}\{u\theta-F(\theta)\}
\]
is the rate function: $K_n/n$ satisfies a large-deviation principle
with speed $n$ and rate $I$ \cite{DemboZeitouni1998}, so that $I(u)$ is
the exponential price of forcing the endpoint onto the ray $k\sim un$.
In thermodynamic language, this is the canonical-to-microcanonical
passage familiar from statistical mechanics: the tilt $\theta$ plays
the role of an external field, $F(\theta)$ is the limiting free energy
of the tilted endpoint ensemble, and $I(u)$ is the corresponding
entropy/rate profile of the macroscopic ratio $K_n/n$ (see
\cite{Touchette2009} for this dictionary).
In the one-dominant-singularity quadratic regimes $B>0$, the
description sharpens to the pointwise ray-scale form
\[
  p_{n,k}
  =
  C(k/n)\,n^{-1/2}\,\exp\{-nI(k/n)\}\,\bigl(1+O(n^{-1})\bigr),
\]
uniformly for $k/n$ in compact subsets of $(0,1)$, with an explicit
Gaussian prefactor $C(\cdot)$ determined by the same singularity map.
The exceptional complex-root case $B=0$ explains why the one-saddle
hypothesis is needed: an antipodal Pearson singularity contributes at
the same exponential scale and produces either an exact span-two
lattice constraint or a persistent even--odd interference correction.
Finally, near the minimum $u_0=F'(0)$ one has
$I(u)=(u-u_0)^{2}/(2F''(0))+O\bigl((u-u_0)^{3}\bigr)$, so the central
Gaussian window is recovered as the quadratic approximation of the
global rate profile---the near-equilibrium expansion of the free
energy.

The novelty is not in the use of Legendre transforms or saddlepoint
approximations as such; these are standard tools. The point is that in
the balanced affine Motzkin class every object entering these methods
is controlled explicitly by one and the same Pearson characteristic
geometry: the PDE, the escape time $\tau(x)$, the limit cumulant
generating function $F$, the rate function $I$, and the finite-$n$
saddlepoint profile.

\medskip
\noindent\textbf{Main contributions.}
The paper contributes the following.

\begin{enumerate}\itemsep0.2em

\item We identify the escape time $t=\tau(x)$ of the Pearson
characteristic flow as the moving singularity of the generating
function and as the central object controlling the terminal-height
distribution.

\item We prove that, in the quadratic regimes,
$\tfrac1n\log\{P_n(e^{\theta})/P_n(1)\}\to
F(\theta)=\log\{\tau(1)/\tau(e^{\theta})\}$, locally uniformly and with
rate $O(n^{-1})$: the escape time of the Pearson flow becomes the limit
cumulant generating function of the endpoint law.

\item We prove the $n$-speed large-deviation principle for $K_n/n$ with
strictly convex good rate function $I=F^{*}$, and give the parametric
Pearson representation $u(x)=x/(Q(x)\tau(x))$,
$I(u(x))=u(x)\log x-\log\{\tau(1)/\tau(x)\}$, valid in all
non-degenerate quadratic regimes.

\item We prove a uniform finite-$n$ Daniels saddlepoint approximation
for the point probabilities $p_{n,k}$, with relative error $O(n^{-1})$
on compact subintervals of the interior of the support in the
one-dominant-singularity regimes $B>0$, and deduce the ray-scale form
with explicit Gaussian prefactor displayed above.

\item We describe the exceptional complex-root case $B=0$, where an
antipodal Pearson singularity produces either an exact span-two lattice
constraint or a persistent even--odd interference correction.

\item As structural input, we show that the balanced condition
$\beta_0=C$ is precisely the locality condition for the
generating-function PDE of the affine model, and we solve the balanced
PDE in closed form in all five Pearson regimes.
\end{enumerate}

\medskip
\noindent\textbf{Relation to earlier work.}
The affine Motzkin and Dyck triangles, including several closed-form
enumerative specialisations, were introduced in \cite{Meshkov2010}; the
present paper develops their probabilistic asymptotic theory. For
constant weights, Banderier and Flajolet \cite{Banderier2002} built a
complete analytic theory in the algebraic setting; our results may be
read as an affine-weight counterpart of that programme, with
exponential generating functions and Pearson characteristics replacing
kernel methods and algebraic singularities. On the probabilistic side
we use classical saddlepoint methods
\cite{LugannaniRice1980,BarndorffCox1989,Butler2007,Kolassa2006}. From
the standpoint of analytic combinatorics, the cumulant expansion
$\kappa_n(\theta)=n\psi(\theta)+O(1)$ underlying our analysis is a
quasi-power, and central and local limit theorems in this framework go
back to Hwang \cite{Hwang1996,Hwang1998}; the contribution of the
present setting is that the limit function is explicit,
$\psi(\theta)=-\log\tau(e^{\theta})$, through the Pearson escape time,
and that the saddlepoint analysis can be carried out uniformly with an
$O(n^{-1})$ error across the interior of the support.

\medskip
\noindent\textbf{Organisation of the paper.}
Section~\ref{sec:gf} defines the affine Motzkin model, derives the
balanced local PDE, and records the Pearson closed forms.
Section~\ref{sec:local} develops coefficient asymptotics and shows why
the local Gaussian window does not capture the global profile.
Section~\ref{sec:daniels} introduces tilting and proves the finite-\(n\)
Daniels saddlepoint approximation. Section~\ref{sec:limit} derives the
limit cumulant generating function, the Legendre rate function, the
large-deviation principle, and the ray-scale form with Gaussian
prefactor. Section~\ref{sec:numerics} presents numerical experiments.
Appendix~\ref{app:technical} collects the technical estimates for the
quadratic regimes; Appendices~\ref{app:linear-constant} and
\ref{app:quadratic} record the degenerate regimes and the explicit
quadratic formulae.

\section{Affine Motzkin paths and the Pearson characteristic flow}
\label{sec:gf}

We begin with the weighted Motzkin model. For each \(n\) and \(k\), let
\(w_{n,k}\) denote the total weight of Motzkin paths of length \(n\)
that terminate at height \(k\). The step-by-step structure gives the
three-term recurrence
\begin{equation}
\label{eq:rec}
  w_{n+1,k}
  =
  \alpha_{k-1}w_{n,k-1}
  +\gamma_k w_{n,k}
  +\beta_k w_{n,k+1},
  \qquad 0\le k\le n+1,
\end{equation}
where terms with indices outside the admissible range are omitted. The
boundary conditions are
\begin{equation}
\label{eq:bound_cond}
  w_{0,0}=1,\qquad
  w_{n,k}=0\quad\text{for }k<0\text{ or }k>n.
\end{equation}
Here \(\alpha_k\) and \(\gamma_k\) weight the up-step and level-step
leaving height \(k\), respectively, while \(\beta_k\) weights the
down-step arriving at height \(k\), equivalently the down-step from
height \(k+1\) to height \(k\).

We work with affine height-dependent weights
\begin{equation}
\label{eq:weights}
  \alpha_k=a\,k+\alpha_0,\qquad
  \beta_k=b\,k+\beta_0,\qquad
  \gamma_k=c\,k+\gamma_0,
\end{equation}
where the parameters are non-negative. As in the introduction, we write
\[
  A=a,\qquad B=c,\qquad C=b,\qquad
  Q(x)=Ax^2+Bx+C,\qquad \Delta=B^2-4AC.
\]

Package the rows into polynomials
\[
  P_n(x)=\sum_{k=0}^n w_{n,k}x^k
\]
and into the exponential generating function
\begin{equation}
\label{eq:gf}
  w(x,t)=\sum_{n\ge0}P_n(x)\frac{t^n}{n!}.
\end{equation}
A coefficient extraction from \eqref{eq:rec} gives
\begin{equation}
\label{Peq}
  P_{n+1}(x)
  =
  Q(x)P_n'(x)
  +(\alpha_0x+\gamma_0)P_n(x)
  +(\beta_0-C)\frac{P_n(x)-P_n(0)}{x}.
\end{equation}
The quotient is a polynomial, with value \(P_n'(0)\) at \(x=0\). The
last term is the only nonlocal term in the polynomial recurrence. It
disappears exactly when
\[
  \beta_0=C.
\]
We call this the \emph{balanced} condition.

Passing to exponential generating functions gives
\begin{equation}
\label{eq:PDE}
  \frac{\partial w}{\partial t}
  =
  Q(x)\frac{\partial w}{\partial x}
  +(\alpha_0x+\gamma_0)w
  +(\beta_0-C)\frac{w(x,t)-w(0,t)}{x},
  \qquad w(x,0)=1,
\end{equation}
where the quotient is understood by analytic continuation at \(x=0\).
In the balanced case this becomes the local first-order PDE
\begin{equation}
\label{eq:PDE_balanced}
  \frac{\partial w}{\partial t}
  =
  Q(x)\frac{\partial w}{\partial x}
  +(\alpha_0x+\gamma_0)w,
  \qquad w(x,0)=1.
\end{equation}

The characteristic equations associated with \eqref{eq:PDE_balanced}
are
\[
  \frac{dx}{dt}=-Q(x),
  \qquad
  \frac{d}{dt}\log w=\alpha_0x+\gamma_0.
\]
Eliminating \(t\) gives
\begin{equation}
\label{eq:Pearson}
  \frac{d}{dx}\log w
  =
  -\,\frac{\alpha_0x+\gamma_0}{Q(x)}
  =
  -\,\frac{\alpha_0x+\gamma_0}{Ax^2+Bx+C}.
\end{equation}
This is the Pearson-type logarithmic derivative which drives the paper:
a linear numerator divided by a quadratic denominator. The sign of the
discriminant
\[
  \Delta=B^2-4AC
\]
separates the three quadratic Pearson regimes, while the case \(A=0\)
gives the constant and linear degeneracies.

The balanced PDE can be solved explicitly in all five regimes. We record
the resulting closed forms because their singularities are the starting
point for the asymptotic analysis.

\begin{theor_eng}[Closed forms in the balanced affine case]
\label{thm:closed-forms}
Assume \(\beta_0=C\). Then \(w(x,t)\) is the unique solution, analytic
at \(t=0\), of \eqref{eq:PDE_balanced}. It is given as follows.

\medskip
\noindent\emph{Constant drift: \(A=0,\ B=0\).}
\begin{equation}
\label{eq:A0Beq0}
  w(x,t)
  =
  \exp\!\Big(\alpha_0xt+\frac{\alpha_0C}{2}t^2+\gamma_0t\Big).
\end{equation}

\medskip
\noindent\emph{Linear drift: \(A=0,\ B>0\).}
\begin{equation}
\label{eq:A0Bneq0}
  w(x,t)
  =
  \exp\!\left(
    \frac{\alpha_0}{B}(e^{Bt}-1)
    \left(x+\frac{C}{B}\right)
    +
    \left(\gamma_0-\frac{\alpha_0C}{B}\right)t
  \right).
\end{equation}
In both cases \(A=0\), the function \(w(\cdot,t)\) is entire in \(x\).

\medskip
\noindent\emph{Quadratic drift: \(A>0,\ \Delta>0\).}
Let
\[
  r_1=\frac{-B-\sqrt{\Delta}}{2A},
  \qquad
  r_2=\frac{-B+\sqrt{\Delta}}{2A},
  \qquad r_1<r_2.
\]
Then
\begin{equation}
\label{eq:w-Delta>0-compact}
  w(x,t)
  =
  \exp\!\big[(\alpha_0r_1+\gamma_0)t\big]\,
  \left[
    \frac{r_1-r_2}
    {(x-r_2)-(x-r_1)\exp\!\big(A(r_1-r_2)t\big)}
  \right]^{\alpha_0/A}.
\end{equation}

\medskip
\noindent\emph{Quadratic drift: \(A>0,\ \Delta=0\).}
Let
\[
  r=-\frac{B}{2A},
  \qquad Q(x)=A(x-r)^2.
\]
Then
\begin{equation}
\label{eq:w-Delta0-closed}
  w(x,t)
  =
  \exp\!\big[(\alpha_0r+\gamma_0)t\big]\,
  \big(1-At(x-r)\big)^{-\alpha_0/A}.
\end{equation}

\medskip
\noindent\emph{Quadratic drift: \(A>0,\ \Delta<0\).}
Set
\[
  p=-\frac{B}{2A},
  \qquad
  q=\frac{\sqrt{-\Delta}}{2A}>0,
  \qquad
  Q(x)=A\big((x-p)^2+q^2\big).
\]
With the principal branch \(\arctan\in(-\pi/2,\pi/2)\),
\begin{equation}
\label{eq:w-Delta<0-cos-final}
  w(x,t)
  =
  \exp\!\big[(\alpha_0p+\gamma_0)t\big]\,
  \left[
    \frac{q}
    {\sqrt{(x-p)^2+q^2}\,
     \cos\!\Big(Aqt+\arctan\frac{x-p}{q}\Big)}
  \right]^{\alpha_0/A}.
\end{equation}

In the formulas with fractional powers, the branch is chosen so that the
bracketed factor equals \(1\) at \(t=0\).
\end{theor_eng}

\begin{proof}
Coefficient comparison in \eqref{eq:PDE_balanced} shows that any
solution analytic at \(t=0\) has Taylor coefficients satisfying
\eqref{Peq} with \(\beta_0=C\). Since this recurrence determines
\(P_{n+1}\) uniquely from \(P_n\), the analytic solution is unique.

The displayed formulae are obtained by the elementary method of
characteristics, by integrating \(dx/dt=-Q(x)\) and
\(d(\log w)/dt=\alpha_0x+\gamma_0\) in the five cases. Equivalently, one
verifies them directly by substitution into \eqref{eq:PDE_balanced} and
by checking the initial condition \(w(x,0)=1\).
\end{proof}

The role of Theorem~\ref{thm:closed-forms} is not merely to provide
closed forms. In the quadratic regimes \(A>0\), the formulas exhibit a
moving \(t\)-singularity. This singularity is the geometric object that
later becomes the limit cumulant generating function after exponential
tilting.

\begin{lemma}[The singular time]
\label{lem:tau-integral}
Assume \(A>0\). Let \(D\) be the real component from which the auxiliary
flow \(dX/dt=Q(X)\) reaches \(+\infty\):
\[
  D=
  \begin{cases}
    (r_2,\infty), & \Delta>0,\\
    (r,\infty),  & \Delta=0,\\
    \mathbb R,   & \Delta<0.
  \end{cases}
\]
For \(x\in D\), define
\begin{equation}
\label{eq:tau-integral}
  \tau(x)=\int_x^\infty\frac{dy}{Q(y)}.
\end{equation}
Then
\begin{equation}
\label{eq:tau-def}
  \tau(x)=
  \begin{cases}
    \displaystyle
    \frac{1}{A(r_2-r_1)}
    \log\frac{x-r_1}{x-r_2},
      & \Delta>0,\quad x>r_2,\\[1.4ex]
    \displaystyle
    \frac{1}{A(x-r)},
      & \Delta=0,\quad x>r,\\[1.4ex]
    \displaystyle
    \frac{\frac{\pi}{2}-\arctan\frac{x-p}{q}}{Aq},
      & \Delta<0.
  \end{cases}
\end{equation}
If \(\alpha_0>0\), then for each \(x\in D\) the same value
\(\tau(x)\) is the first positive real singular time of the corresponding
quadratic closed form in Theorem~\ref{thm:closed-forms}. If
\(\alpha_0=0\), the function \(\tau\) remains the characteristic escape
time, although the moving singularity in \(w(x,t)\) may be absent.

Moreover,
\begin{equation}
\label{eq:tau-derivative}
  \tau'(x)=-\frac1{Q(x)},\qquad
  \tau(x)>0,\qquad
  \tau(x)\sim\frac{1}{Ax}\quad(x\to+\infty).
\end{equation}
Since \(B,C\ge0\), the component \(D\) contains \((0,\infty)\) in all
three quadratic regimes.

If \(B+C>0\), then
\[
  u(x):=\frac{x}{Q(x)\tau(x)}
\]
satisfies
\begin{equation}
\label{eq:u-monotone}
  u'(x)>0,\qquad x>0,
\end{equation}
and
\[
  \lim_{x\downarrow0}u(x)=0,
  \qquad
  \lim_{x\to\infty}u(x)=1.
\]
In the boundary case \(B=C=0\), one has \(\tau(x)=1/(Ax)\) and
\(u(x)\equiv1\).
\end{lemma}

The proof of Lemma~\ref{lem:tau-integral} is given in
Appendix~\ref{app:technical}. The key point for the rest of the paper is
that, after tilting by \(x=e^\theta\), the function \(u\) becomes
\[
  F'(\theta)=u(e^\theta),
\]
the derivative of the limit cumulant generating function. Thus the
monotonicity and range of \(u\) are precisely the ingredients that make
the Legendre correspondence \(F'(\theta)=u\) globally well behaved.

By contrast, in the degenerate regimes \(A=0\), no moving algebraic
singularity occurs. The function \(w(\cdot,t)\) is entire in \(x\), and
the natural asymptotic scale is sublinear rather than the \(n\)-scale
used in the quadratic large-deviation theory. The details are recorded
in Appendix~\ref{app:linear-constant}.

\medskip
\noindent\textbf{Example: permutations.}
For the specialisation
\[
  \alpha_k=k+1,\qquad
  \beta_k=k+1,\qquad
  \gamma_k=2k+1,
\]
one has
\[
  A=1,\qquad B=2,\qquad C=1,\qquad \Delta=0,\qquad r=-1.
\]
Theorem~\ref{thm:closed-forms} gives
\[
  w(x,t)=\big(1-t(x+1)\big)^{-1},
\]
and therefore
\[
  P_n(x)=n!(1+x)^n,\qquad
  w_{n,k}=n!\binom{n}{k}.
\]
After normalisation,
\[
  K_n\sim\mathrm{Binomial}(n,1/2).
\]
This example will serve as a useful check on the asymptotic formulae
below: the general theory reduces here to the classical binomial
Gaussian window and to the rate function
\[
  I(u)=u\log u+(1-u)\log(1-u)+\log2.
\]

\section{Local asymptotics and the limitation of the Gaussian window}
\label{sec:local}

We now pass from the combinatorial array to a probabilistic viewpoint by
normalising the weights. For each \(n\) set
\begin{equation}
\label{eq:pnk-def}
  S_n:=P_n(1),\qquad
  p_{n,k}=\frac{w_{n,k}}{S_n}
         =\frac{[x^k]P_n(x)}{P_n(1)},\qquad 0\le k\le n,
\end{equation}
and let \(K_n\) denote the terminal height under this probability law:
\[
  \mathbb P\{K_n=k\}=p_{n,k}.
\]
Equivalently, the probability generating function of \(K_n\) is
\[
  \mathbb E[x^{K_n}]=\frac{P_n(x)}{P_n(1)}.
\]
Thus
\begin{equation}
\label{eq:mu-sigma}
  \mu_n=\mathbb E[K_n]=\frac{P_n'(1)}{P_n(1)},\qquad
  \sigma_n^2=\mathrm{Var}(K_n)
  =
  \frac{P_n''(1)}{P_n(1)}+\mu_n-\mu_n^2.
\end{equation}
It will be convenient to use the logarithmic derivative
\[
  D_x:=x\frac{d}{dx}.
\]
Then
\[
  \mu_n=\left. D_x\log P_n(x)\right|_{x=1},
  \qquad
  \sigma_n^2=\left. D_x^2\log P_n(x)\right|_{x=1}.
\]

The closed forms of Theorem~\ref{thm:closed-forms} allow us to extract
\[
  P_n(x)=n![t^n]w(x,t)
\]
by singularity analysis in the \(t\)-plane. In the degenerate regimes
\(A=0\), no moving algebraic singularity is present: the functions
\(w(\cdot,t)\) are entire in \(x\), and the natural scale of the
terminal height is sublinear in \(n\). These cases are recorded in
Appendix~\ref{app:linear-constant}. They do not lead to a non-trivial
\(n\)-speed large-deviation theory.

The quadratic regimes \(A>0\) are different. In these regimes the closed
forms have a moving dominant singularity \(t=\tau(x)\), and this is the
first point at which the Pearson characteristic geometry becomes
asymptotically visible. Throughout this section we assume
\[
  A>0,\qquad \alpha_0>0,\qquad \beta_0=C,
\]
and write
\[
  \nu:=\frac{\alpha_0}{A}.
\]
By Lemma~\ref{lem:tau-integral}, the positive real axis lies in the
component on which \(\tau(x)\) is defined.

The standard coefficient-extraction principle is the following. Cauchy's
formula gives
\[
  P_n(x)=\frac{n!}{2\pi i}
  \oint \frac{w(x,t)}{t^{n+1}}\,dt.
\]
In the quadratic regimes the dominant contribution comes from the
algebraic singularity \(t=\tau(x)\). Thus, rather than a new saddle
calculation in each regime, one obtains a uniform transfer theorem. The
proof is given in Appendix~\ref{app:technical}.

\begin{theor_eng}[Quadratic coefficient asymptotics]
\label{thm:quadratic-transfer}
Assume
\[
  A>0,\qquad \alpha_0>0,\qquad \beta_0=C,
\]
and let \(\nu=\alpha_0/A\). For each \(x>0\), the balanced generating
function has, at its first positive singular time \(t=\tau(x)\), an
exact local factorisation
\begin{equation}
\label{eq:exact-local-factor}
  w(x,t)=h(x,t)
  \left(1-\frac{t}{\tau(x)}\right)^{-\nu},
\end{equation}
where \(h(x,t)\) is analytic and non-zero at \(t=\tau(x)\). Set
\[
  \mathcal H(x):=h(x,\tau(x)).
\]
Then, for every compact set \(K\subset(0,\infty)\),
\begin{equation}
\label{eq:master-Pn}
  P_n(x)=
  \frac{\sqrt{2\pi}}{\Gamma(\nu)}\,
  \mathcal H(x)\,
  n^{\nu-\frac12}
  \left(\frac{n}{e\,\tau(x)}\right)^{n}
  \left(1+O(n^{-1})\right),
\end{equation}
uniformly for \(x\in K\). The estimate also holds uniformly in a
sufficiently small complex neighbourhood of \(K\).

The amplitude is
\begin{equation}
\label{eq:H-amplitudes}
  \mathcal H(x)=
  \begin{cases}
  \displaystyle
  \exp\!\big[(\alpha_0r_1+\gamma_0)\tau(x)\big]\,
  \big(A\tau(x)(x-r_2)\big)^{-\nu},
    & \Delta>0,\\[1.5ex]
  \displaystyle
  \exp\!\big[(\alpha_0r+\gamma_0)\tau(x)\big],
    & \Delta=0,\\[1.5ex]
  \displaystyle
  \exp\!\big[(\alpha_0p+\gamma_0)\tau(x)\big]\,
  \big(A\tau(x)\sqrt{(x-p)^2+q^2}\big)^{-\nu},
    & \Delta<0,
  \end{cases}
\end{equation}
with \(r_1,r_2,r,p,q\) as in Theorem~\ref{thm:closed-forms} and
Lemma~\ref{lem:tau-integral}.
\end{theor_eng}

Taking logarithms in \eqref{eq:master-Pn} gives, locally uniformly in a
complex neighbourhood of \(x=1\),
\[
  \log P_n(x)
  =
  C_n-n\log\tau(x)+\log\mathcal H(x)+O(n^{-1}),
\]
where \(C_n\) is independent of \(x\). Therefore the asymptotic expansion
can be differentiated by Cauchy's estimates. Put
\[
  \chi(x):=-\frac{\tau'(x)}{\tau(x)}
          =
          \frac{1}{Q(x)\tau(x)},
  \qquad
  u(x):=x\chi(x)=\frac{x}{Q(x)\tau(x)}.
\]
Then
\begin{equation}
\label{eq:master-moments}
  \mu_n=n\,\chi(1)+O(1),
\end{equation}
and
\begin{equation}
\label{eq:master-variance}
  \sigma_n^2
  =
  n\,v+O(1),
  \qquad
  v:=u'(1)
    =
    \chi(1)+\chi(1)^2-\frac{\tau''(1)}{\tau(1)}.
\end{equation}
In the non-degenerate quadratic case \(B+C>0\),
Lemma~\ref{lem:tau-integral} gives
\[
  u'(x)>0,\qquad x>0,
\]
and hence
\[
  v=u'(1)>0.
\]
In the boundary case \(B=C=0\), one has \(u\equiv1\), so the leading
variance coefficient vanishes and the macroscopic terminal height
degenerates at \(K_n/n=1\).

\medskip
\noindent\textbf{Check: the permutation specialisation.}
For the permutation row of Table~\ref{tab:classical},
\[
  \alpha_k=k+1,\qquad
  \beta_k=k+1,\qquad
  \gamma_k=2k+1.
\]
Here
\[
  A=1,\qquad B=2,\qquad C=1,\qquad
  \Delta=0,\qquad r=-1,
\]
and therefore
\[
  \tau(x)=\frac{1}{x+1},\qquad
  \mathcal H(x)=1,\qquad
  \nu=1.
\]
Theorem~\ref{thm:quadratic-transfer} gives
\[
  P_n(x)
  =
  \sqrt{2\pi}\,n^{1/2}
  \left(\frac{n(x+1)}{e}\right)^n
  \left(1+O(n^{-1})\right),
\]
which is precisely Stirling's approximation to the exact identity
\[
  P_n(x)=n!(1+x)^n.
\]
Furthermore,
\[
  \chi(1)=\frac12,\qquad
  v=u'(1)=\frac14,
\]
so
\[
  \mu_n=\frac n2+O(1),
  \qquad
  \sigma_n^2=\frac n4+O(1),
\]
in agreement with the exact law
\[
  K_n\sim\mathrm{Binomial}(n,1/2).
\]

The conclusion of this section is deliberately modest. The singularity
\(\tau(x)\) determines the leading growth of \(P_n(x)\) and the first
two moments of the terminal-height distribution. It also identifies the
Gaussian scale near the mean. But this local information is not the
global profile of the distribution.

Indeed, around the mean all non-degenerate regimes look Gaussian to
first order. This local Gaussian window is therefore too coarse to
explain the Pearson structure. The distinctions between the Pearson
regimes appear on ray scales
\[
  k\sim un,
\]
where the movement of \(\tau(x)\) under exponential tilting becomes
essential. To see this global structure, we must replace the fixed
expansion at \(x=1\) by a moving tilted expansion.

This is the purpose of the next section. We introduce the finite-\(n\)
cumulant generating function
\[
  \kappa_n(\theta)=\log P_n(e^\theta),
\]
choose the tilt \(\theta_{n,k}\) so that the tilted mean equals the
target height \(k\), and apply Daniels' lattice saddlepoint method. This
will give a finite-\(n\) approximation that interpolates between the
central Gaussian window and the large-deviation tails. In
Section~\ref{sec:limit}, the same tilted geometry will converge to the
limit cumulant generating function
\[
  F(\theta)=\log\frac{\tau(1)}{\tau(e^\theta)}
\]
and to the Legendre rate function \(I(u)\), which describe the global
endpoint profile.

\section{Finite-$n$ global profiles via Daniels' saddlepoint method}
\label{sec:daniels}

We now implement the global programme outlined at the end of the
previous section. The central object is the finite-\(n\) cumulant
generating function
\[
  \kappa_n(\theta)=\log P_n(e^\theta),\qquad
  \widetilde\kappa_n(\theta)=\kappa_n(\theta)-\kappa_n(0)
  =\log\frac{P_n(e^\theta)}{P_n(1)}.
\]
The moment generating function of \(K_n\) is
\[
  M_n(\theta)=\exp\{\widetilde\kappa_n(\theta)\}.
\]
Tilting by \(\theta\) means replacing \(p_{n,k}\) by
\[
  p_{n,k}^{(\theta)}
  =
  \frac{p_{n,k}e^{\theta k}}{M_n(\theta)}.
\]
Under this tilted law,
\begin{equation}
\label{eq:mean-under-tilt}
  \kappa_n'(\theta)=\mathbb E_\theta[K_n],
  \qquad
  \kappa_n''(\theta)=\mathrm{Var}_\theta(K_n).
\end{equation}
Thus, to estimate \(p_{n,k}\), the natural choice of tilt is the real
solution \(\theta_{n,k}\) of
\begin{equation}
\label{eq:saddle-eq-theta}
  \kappa_n'(\theta_{n,k})=k.
\end{equation}
This is the discrete saddlepoint philosophy of Daniels
\cite{Daniels1954}; see also
\cite{BarndorffCox1989,Butler2007,Kolassa2006}.

For later use set
\[
  \psi(\theta)=-\log\tau(e^\theta),
  \qquad
  F(\theta)=\psi(\theta)-\psi(0)
  =
  \log\frac{\tau(1)}{\tau(e^\theta)}.
\]
By Lemma~\ref{lem:tau-integral},
\[
  F'(\theta)=\psi'(\theta)
  =
  u(e^\theta),
  \qquad
  \psi''(\theta)=e^\theta u'(e^\theta).
\]
In the non-degenerate case \(B+C>0\), \(u'(x)>0\) on \((0,\infty)\);
hence \(\psi\) is strictly convex.

We use exact lattice Fourier inversion. For any real \(\theta\),
\begin{equation}
\label{eq:daniels-contour}
  p_{n,k}
  =
  \frac{e^{-k\theta}}{2\pi P_n(1)}
  \int_{-\pi}^{\pi}
    P_n(e^{\theta+is})e^{-iks}\,ds.
\end{equation}
Equivalently, near the real saddle where an analytic branch of
\(\kappa_n(z)=\log P_n(e^z)\) is chosen,
\[
  p_{n,k}
  =
  \frac{1}{2\pi}
  \int_{-\pi}^{\pi}
  \exp\!\left\{
    \kappa_n(\theta+is)-\kappa_n(0)-k(\theta+is)
  \right\}\,ds.
\]
Thus the contour is the finite lattice segment
\(\theta-i\pi\le z\le\theta+i\pi\).

The phase is
\[
  \Phi_{n,k}(\theta)=\kappa_n(\theta)-k\theta,
\]
and its critical point is exactly \(\theta_{n,k}\). A quadratic expansion
there gives the Daniels point-probability approximation
\begin{equation}
\label{eq:daniels}
  p_{n,k}
  =
  \frac{1}{\sqrt{2\pi\,\kappa_n''(\theta_{n,k})}}\,
  \exp\!\left\{
    \kappa_n(\theta_{n,k})-\kappa_n(0)-k\theta_{n,k}
  \right\}
  \bigl(1+R_{n,k}\bigr).
\end{equation}
The main point of this section is to justify \eqref{eq:daniels}
uniformly on compact subintervals of the interior of the support.

The following separation statement is the analytic form of aperiodicity.
It excludes competing \(t\)-singularities of the same modulus after the
Fourier deformation \(x\mapsto xe^{is}\). The hypothesis \(B>0\) is
essential.

\begin{lemma}[Separation of the dominant singularity]
\label{lem:singularity-separation}
Assume
\[
  A>0,\qquad B>0,\qquad \alpha_0>0.
\]
Let \(I\subset(0,\infty)\) be compact and let \(\delta\in(0,\pi)\).
Then there exists \(\eta>0\) such that, for every
\[
  x\in I,\qquad \delta\le |s|\le\pi,
\]
all \(t\)-singularities of \(w(xe^{is},t)\) have modulus at least
\[
  (1+\eta)\tau(x).
\]
\end{lemma}

\noindent
The proof is given in Appendix~\ref{app:technical}. The idea is to use
coefficientwise majorisation to reduce the problem to excluding equality
of singular moduli. In the real-root regimes this becomes a logarithmic
separation statement; in the complex-root regime it becomes a cotangent
disk inequality. The condition \(B>0\) is exactly what rules out the
antipodal equality.

The separation lemma implies the Fourier decay estimate needed for the
saddlepoint proof.

\begin{lemma}[Fourier decay away from the saddle]
\label{lem:fourier-decay}
Assume
\[
  \beta_0=C,\qquad A>0,\qquad \alpha_0>0,\qquad B>0.
\]
Let \(\Theta\subset\mathbb R\) be compact and let
\(\delta\in(0,\pi)\). Then there exist constants \(c>0\), \(C>0\), and
\(N_0\) such that, for all \(n\ge N_0\),
\begin{equation}
\label{eq:fourier-decay}
  \sup_{\theta\in\Theta}
  \sup_{\delta\le |s|\le\pi}
  \left|
  \frac{P_n(e^{\theta+is})}{P_n(e^\theta)}
  \right|
  \le C e^{-cn}.
\end{equation}
\end{lemma}

\noindent
The proof follows from Lemma~\ref{lem:singularity-separation} together
with the uniform transfer estimate of Theorem~\ref{thm:quadratic-transfer};
see Appendix~\ref{app:technical}.

We can now state the uniform saddlepoint approximation.

\begin{theor_eng}[Uniform finite-$n$ saddlepoint approximation]
\label{thm:daniels-uniform}
Assume
\[
  \beta_0=C,\qquad A>0,\qquad \alpha_0>0,\qquad B>0.
\]
Let \(\varepsilon\in(0,\tfrac12)\). Choose a compact interval
\[
  \Theta_\varepsilon=[\theta_-,\theta_+]
\]
such that
\[
  F'(\theta_-)<\varepsilon,
  \qquad
  F'(\theta_+)>1-\varepsilon .
\]
Then there exist constants \(N_0\) and \(C_\varepsilon>0\) such that, for
all \(n\ge N_0\) and all integers
\[
  k\in[\varepsilon n,(1-\varepsilon)n],
\]
the coefficient \(w_{n,k}\) is positive, the saddlepoint equation
\[
  \kappa_n'(\theta_{n,k})=k
\]
has a unique real solution \(\theta_{n,k}\in\Theta_\varepsilon\), and
\begin{equation}
\label{eq:daniels-uniform}
  p_{n,k}
  =
  \frac{1}{\sqrt{2\pi\,\kappa_n''(\theta_{n,k})}}\,
  \exp\!\left\{
    \kappa_n(\theta_{n,k})-\kappa_n(0)-k\theta_{n,k}
  \right\}
  \left(1+R_{n,k}\right),
  \qquad
  |R_{n,k}|\le \frac{C_\varepsilon}{n}.
\end{equation}
The estimate is uniform for \(k/n\in[\varepsilon,1-\varepsilon]\).
\end{theor_eng}

\noindent
The proof is a three-arc saddlepoint argument based on the finite
Fourier inversion formula \eqref{eq:daniels-contour}. The central arc
gives the Gaussian integral, the intermediate arc is controlled by
uniform convexity, and the distant arc is controlled by
Lemma~\ref{lem:fourier-decay}. The details are given in
Appendix~\ref{app:technical}.

The central Gaussian local limit theorem now follows as the small-tilt
limit of the saddlepoint approximation.

\begin{corollary}[Central Gaussian local limit law]
\label{cor:central-llt}
Assume
\[
  \beta_0=C,\qquad A>0,\qquad \alpha_0>0,\qquad B>0.
\]
Uniformly for
\[
  k=\mu_n+O(\sigma_n),
\]
one has
\begin{equation}
\label{eq:central-llt}
  p_{n,k}
  =
  \frac{1}{\sqrt{2\pi\sigma_n^2}}
  \exp\!\left(
    -\frac{(k-\mu_n)^2}{2\sigma_n^2}
  \right)
  \left(1+O(n^{-1/2})\right).
\end{equation}
If \(k-\mu_n=O(1)\), the relative error improves to \(O(n^{-1})\).
\end{corollary}

\noindent
This is the small-tilt limit of
Theorem~\ref{thm:daniels-uniform}; the short expansion is recorded in
Appendix~\ref{app:technical}.

\begin{remark}[The exceptional case \(B=0\) in the complex-root regime]
\label{rem:B0-interference}
The hypothesis \(B>0\) in Theorem~\ref{thm:daniels-uniform} is essential.
When \(B=0\) and \(A,C>0\), one is in the complex-root regime
\(\Delta<0\) with \(p=0\). Then the singularities associated with
\(x\) and \(-x\) may have the same modulus, and the Fourier ratio
\[
  \frac{P_n(-x)}{P_n(x)}
\]
need not decay exponentially.

If \(\gamma_0=0\), there are no level steps in the alternating
specialisation of Table~\ref{tab:classical}. In that case the terminal
height has the exact parity constraint
\[
  K_n\equiv n \pmod 2.
\]
The one-saddle formula must then be replaced by the span-two lattice
version, with an additional factor \(2\) on the admissible sublattice and
zero probability off it.

If \(\gamma_0>0\), the distribution is aperiodic but still exhibits a
persistent even--odd interference. The leading approximation has a
two-singularity form of the type
\[
  p_{n,k}
  \approx
  p_{n,k}^{\mathrm{Daniels}}
  \left(
    1+(-1)^{\,n-k}
    \exp\{-2\gamma_0\tau(e^{\theta_{n,k}})\}
  \right),
\]
up to lower-order corrections. Thus the failure of the one-saddle
formula in this case reflects a genuine antipodal Pearson singularity,
rather than a limitation of the saddlepoint method itself.
\end{remark}

The Daniels formula \eqref{eq:daniels-uniform} therefore provides a
single finite-\(n\) expression for \(p_{n,k}\) in the
one-dominant-singularity quadratic regimes \(B>0\), valid uniformly on
compact subintervals of the interior of the support. The exceptional
complex-root case \(B=0\) is governed by the interference phenomenon
described in Remark~\ref{rem:B0-interference}.

The same tilted geometry will now be passed to the limit \(n\to\infty\).
The finite-\(n\) saddle equation
\[
  \kappa_n'(\theta_{n,k})=k
\]
will become the limiting equation
\[
  F'(\theta)=u,
\]
where
\[
  F(\theta)=\log\frac{\tau(1)}{\tau(e^\theta)}.
\]
The Legendre transform of \(F\) is the large-deviation rate function for
the ray scale \(k\sim un\).

\section{From the moving singularity to the global rate function}
\label{sec:limit}

We now pass from the finite-\(n\) cumulant generating function
\[
  \widetilde\kappa_n(\theta)
  =
  \log\frac{P_n(e^\theta)}{P_n(1)}
\]
to its \(n\)-scale limit. This limit describes the exponential profile
of the terminal height \(K_n\) on ray scales \(K_n\sim un\).

Throughout this section we work in the balanced quadratic case
\[
  \beta_0=C,\qquad A>0,\qquad \alpha_0>0.
\]
The characteristic escape time \(\tau(x)\) is the function introduced in
Lemma~\ref{lem:tau-integral},
\[
  \tau(x)=\int_x^\infty \frac{dy}{Q(y)},\qquad Q(x)=Ax^2+Bx+C,
\]
on the component containing the positive real axis. By
Lemma~\ref{lem:tau-integral}, this component contains \((0,\infty)\).

Define
\[
  \kappa_n(\theta)=\log P_n(e^\theta),
  \qquad
  F(\theta):=
  \lim_{n\to\infty}\frac1n\widetilde\kappa_n(\theta),
\]
whenever the limit exists. The next lemma identifies this limit.

\begin{lemma}[Cumulant expansion and limit CGF]
\label{lem:limit-CGF}
Assume
\[
  \beta_0=C,\qquad A>0,\qquad \alpha_0>0.
\]
Put
\[
  G(\theta)=\log\frac{\mathcal H(e^\theta)}{\mathcal H(1)}.
\]
Then, locally uniformly for \(\theta\in\mathbb R\),
\begin{equation}
\label{eq:kappa-limit-expansion}
  \widetilde\kappa_n(\theta)
  =
  nF(\theta)+G(\theta)+O(n^{-1}),
\end{equation}
where
\begin{equation}
\label{eq:F-def}
  F(\theta)
  =
  \log\frac{\tau(1)}{\tau(e^\theta)}.
\end{equation}
More precisely, the expansion holds uniformly in a sufficiently small
complex neighbourhood of every compact interval in \(\mathbb R\), with
the analytic branches of \(\tau(e^\theta)\), \(\mathcal H(e^\theta)\),
\(F\), and \(G\) chosen by continuation from the real axis. In particular,
\[
  \frac1n\,\widetilde\kappa_n(\theta)
  =
  F(\theta)+O(n^{-1})
\]
locally uniformly in \(\theta\).
\end{lemma}

\begin{proof}
By Theorem~\ref{thm:quadratic-transfer}, applied uniformly for
\(x=e^\theta\) in a sufficiently small complex neighbourhood of any
compact subinterval of \((0,\infty)\),
\[
  P_n(x)
  =
  \frac{\sqrt{2\pi}}{\Gamma(\nu)}
  \mathcal H(x)\,
  n^{\nu-\frac12}
  \left(\frac{n}{e\,\tau(x)}\right)^n
  \left(1+O(n^{-1})\right).
\]
Taking logarithms gives
\[
  \log P_n(x)
  =
  c_n-n\log\tau(x)+\log\mathcal H(x)+O(n^{-1}),
\]
where \(c_n\) is independent of \(x\). Subtracting the value at \(x=1\)
with \(x=e^\theta\) gives \eqref{eq:kappa-limit-expansion}.
\end{proof}

The geometry of \(F\) is inherited directly from the escape-time
function \(\tau\).

\begin{lemma}[Geometry of \(F\)]
\label{lem:geom-F}
Assume
\[
  A>0,\qquad B+C>0.
\]
Then \(F\in C^\infty(\mathbb R)\), and
\[
  F'(\theta)
  =
  u(e^\theta)
  =
  \frac{e^\theta}{Q(e^\theta)\tau(e^\theta)}.
\]
Moreover,
\[
  F''(\theta)=e^\theta u'(e^\theta)>0,
\]
so \(F\) is strictly convex. Finally,
\begin{equation}
\label{eq:Fprime-range}
  \lim_{\theta\to-\infty}F'(\theta)=0,\qquad
  \lim_{\theta\to+\infty}F'(\theta)=1.
\end{equation}
Thus \(F'\) is a bijection from \(\mathbb R\) onto \((0,1)\).
\end{lemma}

\begin{proof}
By Lemma~\ref{lem:tau-integral},
\[
  \tau'(x)=-\frac1{Q(x)}.
\]
Hence
\[
  F'(\theta)
  =
  -\frac{\tau'(e^\theta)e^\theta}{\tau(e^\theta)}
  =
  \frac{e^\theta}{Q(e^\theta)\tau(e^\theta)}
  =
  u(e^\theta).
\]
Differentiating once more gives
\[
  F''(\theta)=e^\theta u'(e^\theta).
\]
When \(B+C>0\), Lemma~\ref{lem:tau-integral} gives
\(u'(x)>0\) for all \(x>0\). Hence \(F''(\theta)>0\) for all
\(\theta\).

The endpoint limits in \eqref{eq:Fprime-range} are also exactly those
of Lemma~\ref{lem:tau-integral}:
\[
  \lim_{x\downarrow0}u(x)=0,\qquad
  \lim_{x\to\infty}u(x)=1.
\]
Since \(x=e^\theta\), this proves the stated range of \(F'\).
\end{proof}

\begin{remark}[The boundary case \(B=C=0\)]
\label{rem:BC0-degenerate}
If \(B=C=0\), then \(Q(x)=Ax^2\) and
\[
  \tau(x)=\frac{1}{Ax}.
\]
Consequently
\[
  F(\theta)=\theta,\qquad F'(\theta)\equiv1.
\]
Thus the \(n\)-scale terminal height degenerates at \(K_n/n=1\). The
Legendre transform is
\[
  I(u)=
  \begin{cases}
    0, & u=1,\\
    +\infty, & u\ne1.
  \end{cases}
\]
For this reason, the non-degenerate large-deviation theory below is
stated under the assumption \(B+C>0\).
\end{remark}

From a probabilistic point of view, \(F'(\theta)\) describes the limiting
terminal height per unit length under exponential tilting. Indeed, under
the \(\theta\)-tilted law,
\[
  \mathbb E_\theta\!\left[\frac{K_n}{n}\right]
  =
  \frac1n\kappa_n'(\theta)
  \longrightarrow
  F'(\theta).
\]
By Lemma~\ref{lem:geom-F}, as \(\theta\) ranges over \(\mathbb R\), the
values \(F'(\theta)\) cover the whole interval \((0,1)\). Thus every
macroscopic terminal-height ratio \(u\in(0,1)\) is represented by a
unique tilt.

We now pass from the tilt parameter \(\theta\) to the macroscopic
terminal-height ratio
\[
  U_n:=\frac{K_n}{n}.
\]
The passage is made by Legendre duality.

\begin{theor_eng}[Large deviations and the rate function]
\label{thm:LDP}
Assume
\[
  \beta_0=C,\qquad A>0,\qquad \alpha_0>0,\qquad B+C>0.
\]
Then \(U_n=K_n/n\) satisfies a large-deviation principle on \([0,1]\)
with speed \(n\) and good convex rate function
\begin{equation}
\label{eq:I-def}
  I(u)
  =
  \sup_{\theta\in\mathbb R}\{u\theta-F(\theta)\}.
\end{equation}
For \(0<u<1\), this takes the explicit Legendre form
\[
  I(u)
  =
  u\,\theta(u)-F(\theta(u)),
  \qquad
  F'(\theta(u))=u.
\]
The function \(I\) is strictly convex on \((0,1)\), and its unique
minimum is attained at
\[
  u_0=F'(0)=u(1)=\frac{1}{Q(1)\tau(1)}.
\]
Equivalently, since \(u(1)=\chi(1)\), the minimiser is
\[
  u_0=\chi(1).
\]
\end{theor_eng}

\begin{proof}
Lemma~\ref{lem:limit-CGF} gives locally uniform convergence
\[
  \frac1n\widetilde\kappa_n(\theta)\to F(\theta)
\]
for all \(\theta\in\mathbb R\). The function \(F\) is finite and
differentiable on all of \(\mathbb R\). The G\"artner--Ellis theorem
therefore gives the large-deviation upper bound for closed sets and the
lower bound at all exposed points of the Legendre transform
\cite[Thm.~2.3.6]{DemboZeitouni1998}. Since Lemma~\ref{lem:geom-F}
shows that \(F'\) is a bijection from \(\mathbb R\) onto \((0,1)\), every
point \(u\in(0,1)\) is exposed and has the unique representing tilt
\(\theta(u)\).

The variables \(U_n\) are supported in the compact interval \([0,1]\),
so exponential tightness is automatic. The endpoint bounds follow by the
standard Chernoff estimates and by lower semicontinuity of the Legendre
transform, obtained by letting \(u\downarrow0\) or \(u\uparrow1\) from
the interior. Thus the LDP holds on all of \([0,1]\) with good rate
function \(I=F^*\).

Strict convexity of \(I\) on \((0,1)\) follows from strict convexity of
\(F\) and the bijectivity of \(F'\). The unique zero of \(I\) is located
at the un-tilted mean slope,
\[
  u_0=F'(0)=u(1)=\chi(1).
\]
\end{proof}

The LDP gives logarithmic probabilities of sets. Point probabilities
require the finite-\(n\) saddlepoint analysis from
Section~\ref{sec:daniels}. Under the stronger one-dominant-singularity
assumption \(B>0\), Theorem~\ref{thm:daniels-uniform} gives a sharper
ray-scale estimate with an explicit Gaussian prefactor.

\begin{corollary}[Ray-scale form with Gaussian prefactor]
\label{cor:ray-prefactor}
Assume
\[
  \beta_0=C,\qquad A>0,\qquad \alpha_0>0,\qquad B>0.
\]
Let \(J\subset(0,1)\) be compact. Put
\[
  G(\theta)=\log\frac{\mathcal H(e^\theta)}{\mathcal H(1)},
  \qquad
  C(u)=\frac{\exp\{G(\theta(u))\}}{\sqrt{2\pi\,F''(\theta(u))}},
\]
where \(\theta(u)\) is defined by \(F'(\theta(u))=u\). Then, uniformly
for integers \(k\) with \(u=k/n\in J\),
\begin{equation}
\label{eq:ray-prefactor}
  p_{n,k}
  =
  C(u)\,n^{-1/2}\,
  \exp\{-nI(u)\}\,
  \bigl(1+O(n^{-1})\bigr).
\end{equation}
\end{corollary}

\noindent
The proof is obtained by substituting the cumulant expansion
\eqref{eq:kappa-limit-expansion} into the uniform Daniels formula
\eqref{eq:daniels-uniform}; details are given in
Appendix~\ref{app:technical}.

\begin{corollary}[Point probabilities on compact interior rays]
\label{cor:pointwise-LDP}
Assume
\[
  \beta_0=C,\qquad A>0,\qquad \alpha_0>0,\qquad B>0.
\]
Let \(J\subset(0,1)\) be compact. Then, uniformly for integers \(k\)
with \(k/n\in J\),
\begin{equation}
\label{eq:pointwise-LDP}
  \log p_{n,k}
  =
  -n\,I(k/n)+O(\log n).
\end{equation}
Equivalently,
\[
  p_{n,k}
  =
  \exp\{-nI(k/n)+O(\log n)\}.
\]
In particular, if \(k=\lfloor un\rfloor\) with \(u\in J\), then
\[
  p_{n,\lfloor un\rfloor}
  =
  \exp\{-nI(u)+o(n)\},
\]
uniformly in \(u\in J\).
\end{corollary}

\begin{proof}
Take logarithms in Corollary~\ref{cor:ray-prefactor}. Since \(J\) is
compact and \(C(u)\) is continuous and positive on \(J\), the quantity
\(\log C(k/n)\) is uniformly bounded for \(k/n\in J\). Hence
\[
  \log p_{n,k}
  =
  -nI(k/n)-\frac12\log n+O(1),
\]
uniformly for \(k/n\in J\), which gives \eqref{eq:pointwise-LDP}.

For the floor statement, if \(k=\lfloor un\rfloor\), then
\(k/n=u+O(n^{-1})\). Since \(I\) is smooth, hence locally Lipschitz, on
a neighbourhood of \(J\), we have
\[
  I(k/n)=I(u)+O(n^{-1})
\]
uniformly in \(u\in J\). Therefore
\[
  nI(k/n)=nI(u)+O(1),
\]
and the claimed \(o(n)\)-form follows.
\end{proof}

\begin{remark}[The prefactor and the logarithmic estimate]
Corollary~\ref{cor:pointwise-LDP} is the logarithmic form of
\eqref{eq:ray-prefactor}. Taking logarithms gives
\[
  \log p_{n,k}
  =
  -nI(k/n)-\frac12\log n+\log C(k/n)+O(n^{-1}),
\]
uniformly for \(k/n\in J\). Thus the \(O(\log n)\) discrepancy in
\eqref{eq:pointwise-LDP} is precisely the Gaussian prefactor on the
ray scale.
\end{remark}

\begin{remark}[Point probabilities when \(B=0\)]
The pointwise estimates in Corollaries~\ref{cor:ray-prefactor} and
\ref{cor:pointwise-LDP} are stated under \(B>0\), matching
Theorem~\ref{thm:daniels-uniform}. When \(B=0\) in the complex-root
regime, the LDP of Theorem~\ref{thm:LDP} is still governed by the same
Legendre transform, but point probabilities may require the
two-singularity correction described in Remark~\ref{rem:B0-interference}.
In the exactly periodic case \(\gamma_0=0\), probabilities vanish off
the admissible parity sublattice, so no unrestricted pointwise estimate
can hold.
\end{remark}

For explicit computations it is useful to return to the \(x\)-variable.
Write
\[
  x=e^\theta.
\]
Then
\[
  F'(\theta)=u(x)
  =
  \frac{x}{Q(x)\tau(x)}.
\]
Thus the rate function has the parametric representation
\begin{equation}
\label{eq:uI-param}
  u=u(x)=\frac{x}{Q(x)\tau(x)},\qquad
  I(u(x))
  =
  u(x)\log x
  -
  \log\frac{\tau(1)}{\tau(x)}.
\end{equation}
This parametrisation is valid for all \(x>0\) in the non-degenerate
quadratic case \(B+C>0\), and \(x\mapsto u(x)\) maps \((0,\infty)\)
bijectively onto \((0,1)\).

Differentiating the Legendre relation gives
\[
  I'(u)=\theta(u),
  \qquad
  I''(u)=\frac{1}{F''(\theta(u))}.
\]
In particular, at the minimum \(u_0=F'(0)\),
\[
  I''(u_0)=\frac{1}{F''(0)}.
\]
Consequently, near the typical value \(u_0=F'(0)\),
\[
  I(u)
  =
  \frac{(u-u_0)^2}{2F''(0)}
  +O\bigl((u-u_0)^3\bigr).
\]
Thus the central Gaussian window is the quadratic approximation to the
global rate profile. In this sense the local normal law is only the
local shadow of the Pearson-driven large-deviation geometry.
Also,
\[
  \frac{1}{n}\operatorname{Var}(K_n)
  =
  \frac{\kappa_n''(0)}{n}
  \longrightarrow
  F''(0),
\]
in agreement with the central window in
Corollary~\ref{cor:central-llt}.

\medskip
\noindent\textbf{Pearson parametrisations of the rate function.}
The general parametrisation \eqref{eq:uI-param} becomes explicit in each
quadratic Pearson regime.

\medskip
\noindent\emph{Two real roots: \(\Delta>0\).}
Let
\[
  d=r_2-r_1,\qquad
  L(x)=\log\frac{x-r_1}{x-r_2},\qquad
  S(x)=(x-r_1)(x-r_2).
\]
Then
\[
  \tau(x)=\frac{L(x)}{Ad},
  \qquad
  u(x)=\frac{dx}{S(x)L(x)}.
\]
Therefore
\[
  I(u(x))
  =
  u(x)\log x-\log\frac{L(1)}{L(x)}.
\]

\medskip
\noindent\emph{Double root: \(\Delta=0\).}
Let
\[
  Q(x)=A(x-r)^2.
\]
In the non-degenerate case \(B+C>0\), necessarily \(r<0\). Then
\[
  \tau(x)=\frac{1}{A(x-r)},\qquad
  F(\theta)=\log\frac{e^\theta-r}{1-r},
\]
and
\[
  u(x)=\frac{x}{x-r}.
\]
Solving for \(x\) gives
\[
  x=\frac{(-r)u}{1-u},\qquad 0<u<1.
\]
Substitution into the Legendre formula yields
\begin{equation}
\label{eq:Iu-Delta0}
  I(u)
  =
  u\log u
  +(1-u)\log(1-u)
  +(u-1)\log(-r)
  +\log(1-r),
  \qquad 0<u<1.
\end{equation}
For the permutation specialisation \(r=-1\), this reduces to
\[
  I(u)=u\log u+(1-u)\log(1-u)+\log2,
\]
the Kullback--Leibler divergence from the Bernoulli law with parameter
\(1/2\).

If \(r=0\), equivalently \(B=C=0\), one is in the degenerate boundary
case of Remark~\ref{rem:BC0-degenerate}, and
\[
  I(1)=0,\qquad I(u)=+\infty\quad(u\ne1).
\]

\medskip
\noindent\emph{Complex conjugate roots: \(\Delta<0\).}
Let
\[
  \Theta(x)=\frac{\pi}{2}-\arctan\frac{x-p}{q}.
\]
Then
\[
  \tau(x)=\frac{\Theta(x)}{Aq},
  \qquad
  u(x)=
  \frac{xq}{\big((x-p)^2+q^2\big)\Theta(x)}.
\]
Therefore
\[
  I(u(x))
  =
  u(x)\log x-\log\frac{\Theta(1)}{\Theta(x)}.
\]

\medskip
In summary, the entire non-degenerate large-deviation profile is governed
by the Pearson escape time \(\tau(x)\). The Pearson geometry is
logarithmic when \(\Delta>0\), rational when \(\Delta=0\), and
trigonometric when \(\Delta<0\); the same three forms pass from
\(\tau(x)\) to \(F(\theta)\) and then to the global rate function
\(I(u)\).

\section{Numerical computation and results}
\label{sec:numerics}

To illustrate the theory we report numerical results for the finite-\(n\)
terminal-height distribution and compare the three analytic descriptions
developed above: the central Gaussian window, the finite-\(n\)
saddlepoint approximation, and the large-deviation profile obtained from
the limit cumulant generating function \(F\).

Exact probabilities were computed from the recurrence for the polynomials
\(P_n(x)\), followed by normalisation by \(P_n(1)\). The Daniels curve
was computed by solving the saddlepoint equation
\[
  \kappa_n'(\theta_{n,k})=k
\]
and substituting the resulting \(\theta_{n,k}\) into the saddlepoint
formula \eqref{eq:daniels}. Near the endpoints, where the raw
point-saddle equation becomes numerically ill-conditioned, the plotted
curve uses the standard Lugannani--Rice endpoint stabilisation \cite{LugannaniRice1980}. This
endpoint stabilisation is used only for the numerical display; the
uniform theorem proved above applies on compact subintervals of the
interior of the support.

\begin{figure}[t]
  \centering
  \includegraphics[width=0.70\textwidth]{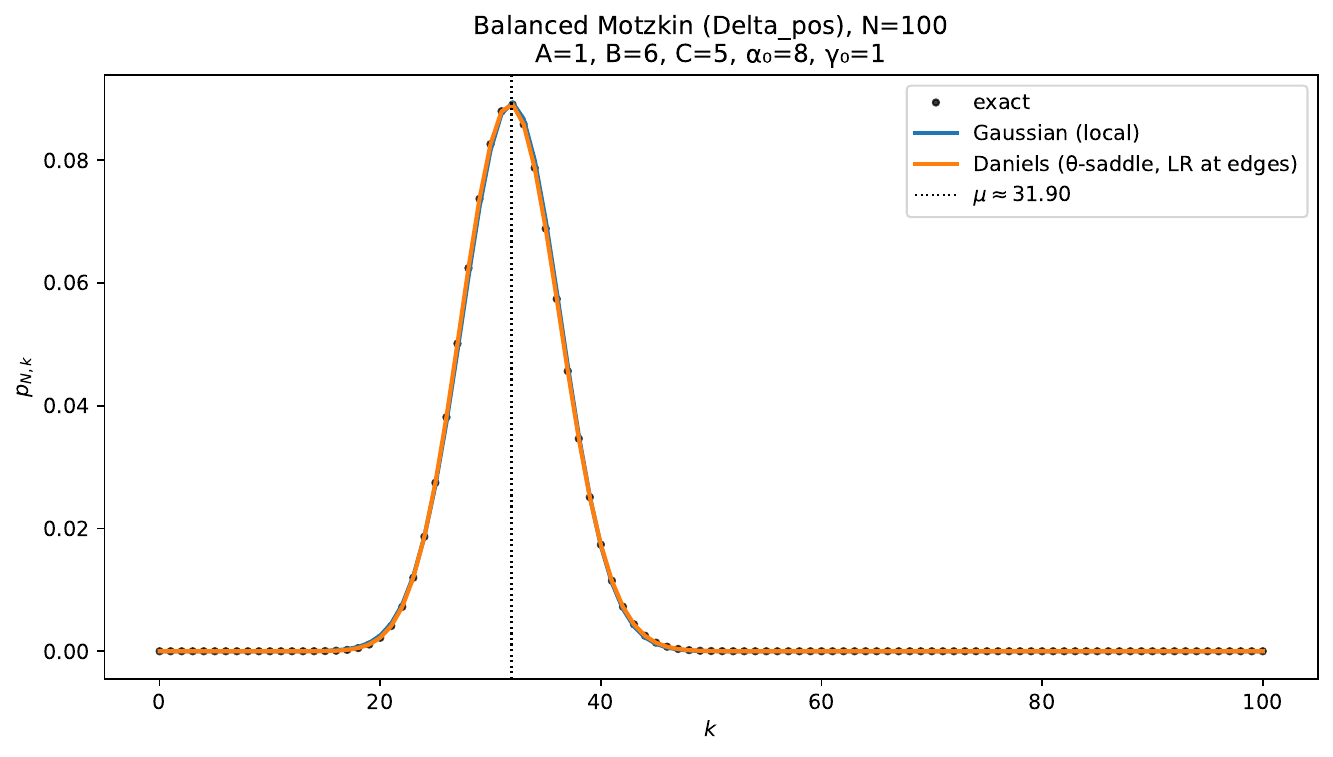}
  \caption{Two real roots (\(\Delta>0\)), linear scale. Exact distribution,
  Gaussian local approximation, and Daniels saddlepoint curve, as indicated in
  the legend. Parameters:
  \(A=1\), \(B=6\), \(C=5\), \(\alpha_0=8\), \(\gamma_0=1\),
  \(\beta_0=C\), and \(n=100\).}
  \label{fig:profile-DeltaPos-lin}
\end{figure}

\begin{figure}[t]
  \centering
  \includegraphics[width=0.70\textwidth]{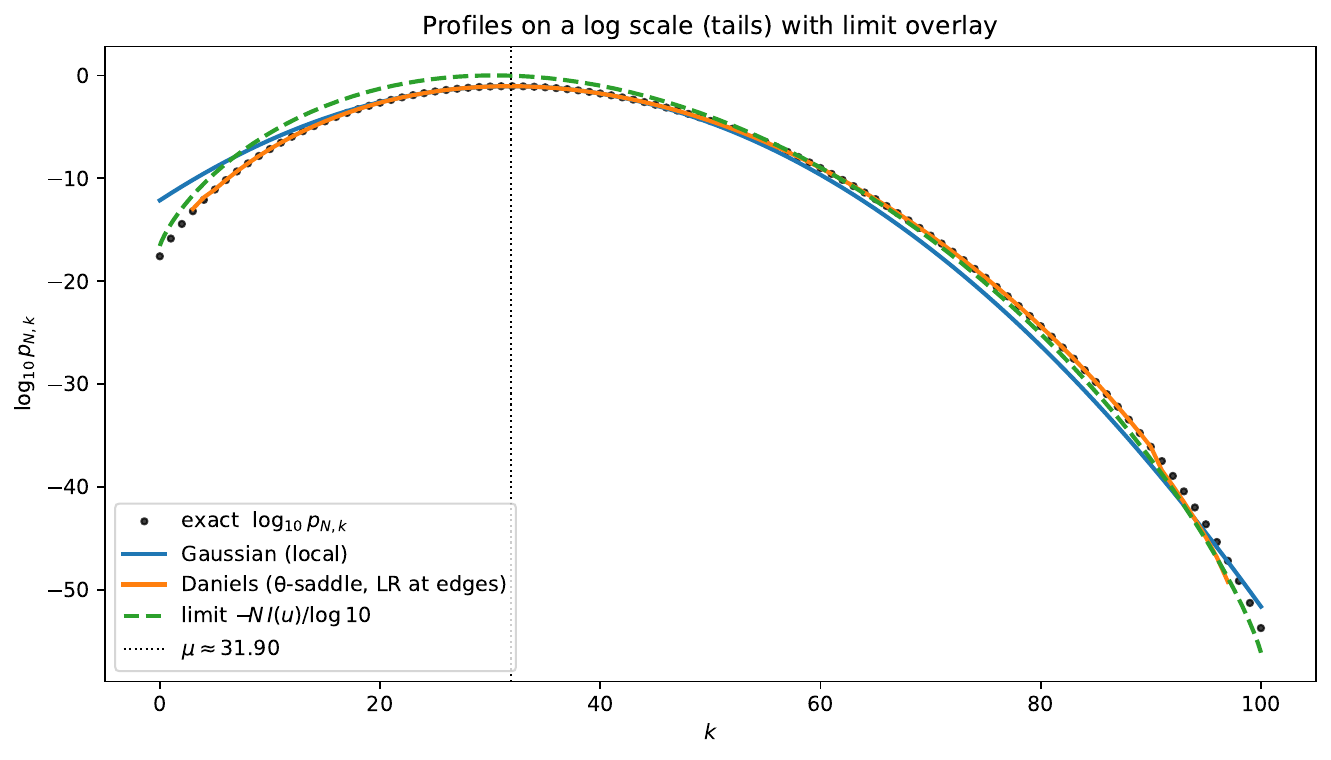}
  \caption{Two real roots (\(\Delta>0\)), logarithmic scale. The same
  distribution as in Figure~\ref{fig:profile-DeltaPos-lin}, with the
  large-deviation overlay \(-nI(k/n)/\log 10\).}
  \label{fig:profile-DeltaPos-log}
\end{figure}

\begin{figure}[t]
  \centering
  \includegraphics[width=0.70\textwidth]{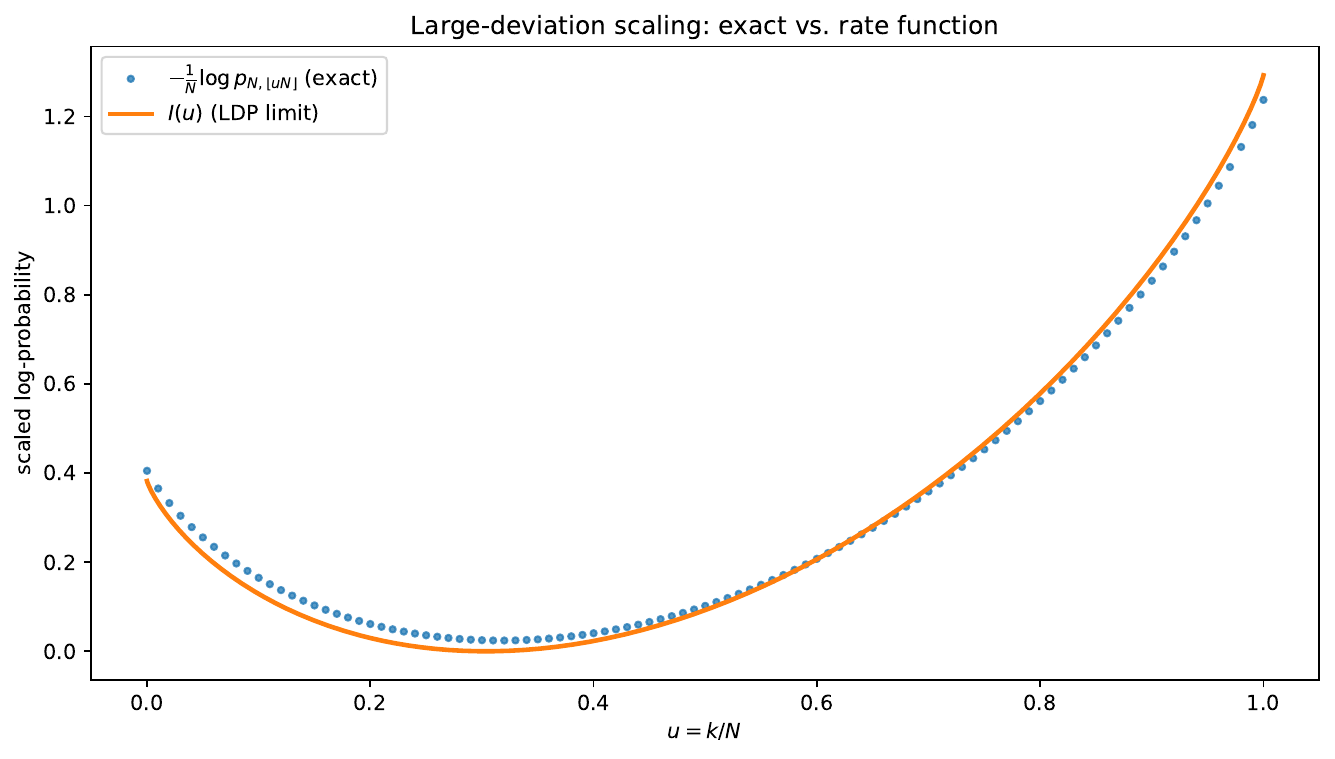}
\caption{Large-deviation scaling for the same \(\Delta>0\) example, with
\(N=100\). Points show \(-N^{-1}\log p_{N,\lfloor uN\rfloor}\), and the
solid curve is the rate function \(I(u)\). The vertical discrepancy is
the Gaussian-prefactor correction predicted by
Corollary~\ref{cor:ray-prefactor}, of order \(O(\log N/N)\).}
  \label{fig:ldp-Delta-pos}
\end{figure}

Figures~\ref{fig:profile-DeltaPos-lin}--\ref{fig:ldp-Delta-pos} show a
representative one-dominant-singularity example in the two-real-root
regime. The parameters are
\[
  A=1,\qquad B=6,\qquad C=5,\qquad
  \alpha_0=8,\qquad \gamma_0=1,\qquad \beta_0=C,
\]
so that
\[
  \Delta=B^2-4AC=16>0.
\]
In particular \(B>0\), and therefore the hypotheses of
Theorem~\ref{thm:daniels-uniform} are satisfied. The path length is
\(n=100\).

On the linear scale, Figure~\ref{fig:profile-DeltaPos-lin} shows that
the Daniels approximation is visually indistinguishable from the exact
distribution across the peak and shoulders. The Gaussian approximation
captures the central window but loses accuracy in the shoulders, as
expected from a purely local central approximation.

On the logarithmic scale, Figure~\ref{fig:profile-DeltaPos-log} shows
the same distribution over several decades of probability. The Daniels
curve continues to track the exact probabilities deep into the tails.
The large-deviation overlay
\[
  -\,\frac{nI(k/n)}{\log 10}
\]
has the correct exponential scale, but differs from the finite-\(n\)
probabilities by the expected prefactor contribution of order
\(O(\log n)\) on the log-probability scale.

Finally, Figure~\ref{fig:ldp-Delta-pos} compares the scaled exact values
\[
  -\frac{1}{N}\log p_{N,\lfloor uN\rfloor}
\]
with the rate function \(I(u)\). The visible vertical discrepancy is not
an error in the rate function: it is the finite-\(N\) Gaussian prefactor.
Indeed, Corollary~\ref{cor:ray-prefactor} gives, uniformly on compact
subintervals of the interior,
\[
  -\frac{1}{N}\log p_{N,\lfloor uN\rfloor}
  =
  I(u)
  +
  \frac{1}{N}\left\{
    \frac12\log\bigl(2\pi N F''(\theta(u))\bigr)
    -G(\theta(u))
  \right\}
  +O(N^{-2}),
\]
up to the harmless replacement of \(u\) by \(\lfloor uN\rfloor/N\). Thus
the discrepancy is of order \(O(\log N/N)\), as seen in the plot. This
makes visible the geometry underlying the theory: the global profile is
controlled by the Pearson singularity map \(x\mapsto\tau(x)\), from
which both \(F\) and \(I\) are obtained.

In summary, the numerical results illustrate the hierarchy proved in the
paper. The Gaussian approximation is the local quadratic shadow of the
rate function near its minimum; the Daniels saddlepoint formula gives the
finite-\(n\) pointwise profile in the interior; and the rate function
\(I(u)\), together with the prefactor in
Corollary~\ref{cor:ray-prefactor}, gives the ray-scale global shape. In
the present \(\Delta>0\), \(B>0\) example all three descriptions are
projections of the same tilted Pearson geometry.

\section*{Conclusion}

The main message of the paper is that, in the balanced affine Motzkin
model, Pearson geometry controls the endpoint distribution globally. The
balanced condition turns the generating-function equation into a local
first-order PDE; its characteristic flow has Pearson form, and the escape
time
\[
  \tau(x)=\int_x^\infty \frac{dy}{Q(y)}
\]
becomes the central analytic object of the theory.

This single function organises the asymptotics at several levels. It
gives the coefficient growth and the Gaussian scale near the mean; after
tilting, it gives the limit cumulant generating function
\[
  F(\theta)=\log\frac{\tau(1)}{\tau(e^\theta)};
\]
and its Legendre transform gives the large-deviation rate function
\(I(u)\). Thus the local Gaussian law is only the quadratic
approximation to a global Pearson-driven rate profile.

The finite-\(n\) Daniels saddlepoint formula is the corresponding
non-asymptotic version of the same geometry. In the
one-dominant-singularity regimes \(B>0\), it gives a uniform pointwise
profile across compact subintervals of the interior and, in the limit,
reduces to the ray-scale form
\[
  p_{n,k}\sim C(k/n)n^{-1/2}e^{-nI(k/n)}.
\]
The exceptional complex-root case \(B=0\) shows that the one-saddle
hypothesis is sharp: an antipodal Pearson singularity may produce either
a span-two lattice constraint or a persistent even--odd interference
term.

The classical specialisations in Table~\ref{tab:classical} show that
this mechanism is not detached from combinatorics: the same discriminant
trichotomy separates permutations, even alternating permutations, and
ordered set partitions. In this sense the paper identifies a route from
affine weighted lattice-path recurrences to global endpoint profiles:
\[
  \text{Pearson characteristic flow}
  \Longrightarrow
  \tau(x)
  \Longrightarrow
  F(\theta)
  \Longrightarrow
  I(u)
  \Longrightarrow
  \text{finite- and large-\(n\) endpoint asymptotics}.
\]

\newpage
\appendix
\numberwithin{equation}{section}
\numberwithin{figure}{section}
\numberwithin{table}{section}

\section{Technical estimates for the quadratic regimes}
\label{app:technical}

This appendix contains the technical estimates used in
Sections~\ref{sec:gf}--\ref{sec:daniels}. We do not repeat the elementary
characteristic derivation of the closed forms. Instead we prove the
properties of the singularity map, the uniform transfer estimate, and the
Fourier decay estimates needed for the finite-\(n\) saddlepoint theorem.

\subsection{The singularity map and monotonicity}

\begin{proof}[Proof of Lemma~\ref{lem:tau-integral}]
The integral
\[
  \tau(x)=\int_x^\infty\frac{dy}{Q(y)}
\]
is finite for \(A>0\), since \(Q(y)\sim Ay^2\) as \(y\to+\infty\).
Evaluating this integral gives the three displayed formulae in
\eqref{eq:tau-def}. Differentiating the integral representation gives
\[
  \tau'(x)=-\frac1{Q(x)},
\]
and the asymptotic
\[
  \tau(x)\sim\frac1{Ax}
  \qquad (x\to+\infty)
\]
follows from \(Q(x)\sim Ax^2\).

If \(\alpha_0>0\), the same values are the first positive real
singular times in the closed forms of Theorem~\ref{thm:closed-forms}.
For \(\Delta>0\), the denominator in
\eqref{eq:w-Delta>0-compact} vanishes at
\[
  t=
  \frac{1}{A(r_2-r_1)}
  \log\frac{x-r_1}{x-r_2}.
\]
For \(x>r_2\), this is the unique positive real zero. For
\(\Delta=0\), the factor \(1-At(x-r)\) vanishes at
\[
  t=\frac1{A(x-r)}.
\]
For \(\Delta<0\), the first positive zero of the cosine in
\eqref{eq:w-Delta<0-cos-final} occurs when
\[
  Aqt+\arctan\frac{x-p}{q}=\frac{\pi}{2},
\]
which gives the third formula in \eqref{eq:tau-def}.

Since \(B,C\ge0\), the component \(D\) contains \((0,\infty)\). Indeed,
in the two-real-root case,
\[
  r_2=\frac{-B+\sqrt{B^2-4AC}}{2A}\le0,
\]
and in the double-root case \(r=-B/(2A)\le0\).

It remains to prove the monotonicity of
\[
  u(x)=\frac{x}{Q(x)\tau(x)}.
\]
Write
\[
  \frac1{u(x)}
  =
  \frac{Q(x)\tau(x)}{x}
  =
  \int_1^\infty \frac{Q(x)}{Q(xs)}\,ds.
\]
For
\[
  h(z)=\frac{zQ'(z)}{Q(z)}
\]
one computes
\[
  h'(z)
  =
  \frac{ABz^2+4ACz+BC}{Q(z)^2}.
\]
Thus \(h\) is strictly increasing on \((0,\infty)\) whenever
\(B+C>0\).

We may differentiate under the integral sign locally uniformly in
\(x>0\). On each compact interval \(K\subset(0,\infty)\),
\[
  0\le \frac{Q(x)}{Q(xs)}\le \frac{C_K}{s^2},
  \qquad s\ge1,
\]
and
\[
  \left|
  \frac{\partial}{\partial x}
  \frac{Q(x)}{Q(xs)}
  \right|
  \le \frac{C_K'}{s^2},
\]
which is integrable on \([1,\infty)\). Hence
\[
  \left(\frac1{u(x)}\right)'
  =
  \int_1^\infty
  \frac{Q(x)}{Q(xs)}
  \frac{h(x)-h(xs)}{x}\,ds.
\]
For \(s>1\) we have \(xs>x\), and therefore \(h(x)-h(xs)<0\). Hence
\[
  \left(\frac1{u(x)}\right)'<0,
  \qquad
  u'(x)>0.
\]

The endpoint \(u(x)\to1\) as \(x\to\infty\) follows by dominated
convergence:
\[
  \frac{Q(x)}{Q(xs)}\to\frac1{s^2},
  \qquad
  \frac1{u(x)}\to\int_1^\infty\frac{ds}{s^2}=1.
\]
As \(x\downarrow0\), either \(C>0\), in which case
\(Q(x)\tau(x)/x\to\infty\), or \(C=0\) and \(B>0\), in which case
\(\tau(x)\to\infty\) and again \(Q(x)\tau(x)/x\to\infty\). Thus
\(u(x)\to0\). If \(B=C=0\), then \(Q(x)=Ax^2\),
\(\tau(x)=1/(Ax)\), and \(u(x)\equiv1\).
\end{proof}

\subsection{Proof of the quadratic transfer estimate}

\begin{proof}[Proof of Theorem~\ref{thm:quadratic-transfer}]
Let
\[
  \nu=\frac{\alpha_0}{A}.
\]
The closed forms in Theorem~\ref{thm:closed-forms} have, at
\(t=\tau(x)\), an exact local factorisation
\[
  w(x,t)=h(x,t)\left(1-\frac{t}{\tau(x)}\right)^{-\nu},
\]
with \(h\) analytic and non-zero at \(t=\tau(x)\).

We compute the amplitudes \(h(x,\tau(x))\). In the case
\(\Delta>0\), put
\[
  \lambda=A(r_2-r_1)>0.
\]
The denominator in \eqref{eq:w-Delta>0-compact} is
\[
  (x-r_2)-(x-r_1)e^{-\lambda t}.
\]
At \(t=\tau(x)\) it has a simple zero, and
\[
  (x-r_2)-(x-r_1)e^{-\lambda t}
  =
  -\lambda(x-r_2)\tau(x)
  \left(1-\frac{t}{\tau(x)}\right)(1+O(t-\tau(x))).
\]
Thus
\[
  \mathcal H(x)
  =
  \exp\!\big[(\alpha_0r_1+\gamma_0)\tau(x)\big]\,
  \big(A\tau(x)(x-r_2)\big)^{-\nu}.
\]

In the case \(\Delta=0\), the factorisation follows immediately from
\[
  1-At(x-r)=1-\frac{t}{\tau(x)},
\]
and
\[
  \mathcal H(x)
  =
  \exp\!\big[(\alpha_0r+\gamma_0)\tau(x)\big].
\]

In the case \(\Delta<0\), let
\[
  R(x)=\sqrt{(x-p)^2+q^2},
  \qquad
  \phi(x)=\arctan\frac{x-p}{q}.
\]
Since
\[
  Aq\tau(x)+\phi(x)=\frac{\pi}{2},
\]
we have, near \(t=\tau(x)\),
\[
  \cos\!\big(Aqt+\phi(x)\big)
  =
  \sin\!\big(Aq(\tau(x)-t)\big)
  =
  Aq\,\tau(x)
  \left(1-\frac{t}{\tau(x)}\right)(1+O(t-\tau(x))).
\]
Hence
\[
  \mathcal H(x)
  =
  \exp\!\big[(\alpha_0p+\gamma_0)\tau(x)\big]\,
  \big(A\tau(x)R(x)\big)^{-\nu}.
\]

For \(x>0\), the singularity \(t=\tau(x)\) is the unique singularity of
minimal modulus. For \(\Delta=0\) this is immediate. For \(\Delta>0\),
the other singularities are
\[
  \tau(x)+\frac{2\pi i m}{A(r_2-r_1)},
  \qquad m\ne0,
\]
and have larger modulus. For \(\Delta<0\), the zeros of the cosine are
\[
  \tau(x)+\frac{\pi m}{Aq},
  \qquad m\in\mathbb Z.
\]
Since \(x>0\) and \(p\le0\), one has
\[
  0<\tau(x)<\frac{\pi}{2Aq};
\]
therefore the nearest negative zero has modulus larger than \(\tau(x)\),
and all other zeros are farther away.

This separation is uniform for \(x\) in compact subsets of
\((0,\infty)\), and also in a sufficiently small complex neighbourhood
of such compact subsets. Hence a \(\Delta\)-domain at the dominant
singularity \(t=\tau(x)\) can be chosen locally uniformly in \(x\). The
algebraic transfer theorem \cite[Thm.~VI.4]{Flajolet2009}, equivalently
the Flajolet--Odlyzko transfer theorem \cite{FlajoletOdlyzko1990}, gives
\[
  [t^n]w(x,t)
  =
  \frac{\mathcal H(x)}{\Gamma(\nu)}
  \tau(x)^{-n}n^{\nu-1}
  \left(1+O(n^{-1})\right),
\]
locally uniformly in \(x\). Multiplying by Stirling's formula
\[
  n!=\sqrt{2\pi}\,n^{n+\frac12}e^{-n}
  \left(1+O(n^{-1})\right)
\]
gives \eqref{eq:master-Pn}.
\end{proof}

\subsection{Fourier separation and decay}

\begin{proof}[Proof of Lemma~\ref{lem:singularity-separation}]
Put \(z=xe^{is}\). Since the coefficients of \(w(x,t)\) are
non-negative,
\[
  |w(z,t)|\le w(x,|t|)
\]
in the common domain of convergence. Hence the \(t\)-radius of
convergence at \(z\) is at least \(\tau(x)\). It remains to exclude
equality for \(\delta\le |s|\le\pi\).

In the double-root case \(\Delta=0\), the singularity is
\[
  t(z)=\frac1{A(z-r)},\qquad r=-\frac{B}{2A}<0.
\]
For \(s\ne0\),
\[
  |xe^{is}-r|^2
  =
  x^2+r^2-2xr\cos s
  <
  x^2+r^2-2xr
  =
  (x-r)^2,
\]
and therefore \(|t(z)|>\tau(x)\).

In the two-real-root case \(\Delta>0\), write
\[
  r_1<r_2\le0,\qquad
  \lambda=A(r_2-r_1)>0,
\]
and
\[
  R(z)=\frac{z-r_1}{z-r_2}.
\]
The singularities are
\[
  t_m(z)
  =
  \frac{1}{\lambda}
  \bigl(\operatorname{Log}R(z)+2\pi i m\bigr),
  \qquad m\in\mathbb Z,
\]
with the branch of the logarithm obtained by continuation from the
positive real axis. It is enough to consider the branch with minimal
absolute imaginary part.

Put
\[
  u=-r_1,\qquad v=-r_2,\qquad u>v\ge0,
\]
and define, for \(\Im z>0\),
\[
  G(z)=\operatorname{Log}\frac{z+u}{z+v},
  \qquad
  g(z)=\frac{z}{(z+u)(z+v)G(z)}.
\]
The branch is chosen so that \(G(x)>0\) for \(x>0\). The Möbius map
\((z+u)/(z+v)\) sends the upper half-plane into the lower half-plane,
so \(G\) is analytic in \(\Im z>0\) with
\[
  -\pi<\Im G(z)<0.
\]
On the real boundary outside \((-u,-v)\), \(G\) is real and
\(\Im g=0\). On the interval \((-u,-v)\), approached from above,
\[
  G(t)=\log\left|\frac{t+u}{t+v}\right|-i\pi,
\]
and
\[
  \frac{t}{(t+u)(t+v)}>0.
\]
Thus
\[
  \Im g(t)
  =
  \frac{t}{(t+u)(t+v)}
  \frac{\pi}
  {\log^2\left|\frac{t+u}{t+v}\right|+\pi^2}
  >0.
\]
Moreover,
\[
  g(z)\to\frac{1}{u-v}
  \qquad (z\to\infty),
\]
so \(\Im g(z)\to0\) at infinity. At the endpoints
\(z=-u\) and \(z=-v\), the boundary values of \(\Im g\) on the cut
\((-u,-v)\) tend to \(+\infty\). Hence the minimum principle may be
applied on the upper half-plane truncated by a large semicircle and by
small semicircles around the two endpoints, and then the truncation
parameters may be removed. We obtain
\[
  \Im g(z)>0,\qquad \Im z>0.
\]

Now let \(z=xe^{is}\), \(0<s<\pi\). Then
\[
  \frac{d}{ds}G(xe^{is})
  =
  -\,\frac{i xe^{is}(u-v)}
  {(xe^{is}+u)(xe^{is}+v)}.
\]
Consequently,
\[
  \frac{d}{ds}\log|G(xe^{is})|
  =
  (u-v)\Im g(xe^{is})>0.
\]
Hence
\[
  \left|
  \operatorname{Log}\frac{xe^{is}-r_1}{xe^{is}-r_2}
  \right|
  >
  \log\frac{x-r_1}{x-r_2}
  \qquad (s\ne0),
\]
and all singularities have modulus strictly larger than \(\tau(x)\).

In the complex-root case \(\Delta<0\), put
\[
  p=-\frac{B}{2A}<0,\qquad
  q=\frac{\sqrt{-\Delta}}{2A}>0,\qquad
  b=-p>0.
\]
The singularities are the zeros of
\[
  \cos\!\left(Aqt+\arctan\frac{z-p}{q}\right).
\]
Equivalently, with \(T=Aqt\), a singularity satisfies
\[
  z=q\cot T-b.
\]
For real \(x>0\), the positive singularity corresponds to
\[
  T=\delta_x:=Aq\,\tau(x)=\arctan\frac{q}{x+b},
\]
so that
\[
  x=q\cot\delta_x-b,
  \qquad 0<\delta_x<\frac{\pi}{2}.
\]

We claim that no singularity with \(z=xe^{is}\), \(s\ne0\), can have
\(|t|\le\tau(x)\). Suppose otherwise. Then for some \(T\) with
\(|T|\le\delta_x\),
\[
  xe^{is}=q\cot T-b.
\]
We use the elementary inequality
\begin{equation}
\label{eq:cot-disk-ineq-app}
  |\cot T-\beta|
  \ge
  \cot\delta-\beta,
  \qquad |T|\le\delta,
\end{equation}
valid for \(\beta>0\), \(0<\delta<\pi/2\), and
\(\beta\tan\delta<1\), with equality only at \(T=\delta\).

Indeed, set
\[
  \Phi(T)=\frac{1}{\cot T-\beta}
  =
  \frac{\tan T}{1-\beta\tan T}.
\]
Since \(\tan T\) has a power series with non-negative coefficients and
\(\beta\tan\delta<1\), the expansion
\[
  \Phi(T)=\sum_{j\ge0}\beta^j(\tan T)^{j+1}
\]
is absolutely convergent for \(|T|\le\delta\), and gives
\[
  |\Phi(T)|\le \Phi(|T|)\le \Phi(\delta).
\]
This proves \eqref{eq:cot-disk-ineq-app}. Applying it with
\[
  \beta=\frac{b}{q},\qquad \delta=\delta_x,
\]
we get
\[
  |xe^{is}|
  =
  q\left|\cot T-\frac{b}{q}\right|
  \ge
  q\cot\delta_x-b
  =
  x.
\]
Since equality already holds in \(|xe^{is}|=x\), equality must hold in
the inequality; hence \(T=\delta_x\), and therefore \(xe^{is}=x\), so
\(s=0\), a contradiction.

Thus in all three regimes the closest singularity for \(z=xe^{is}\),
\(s\ne0\), has modulus strictly larger than \(\tau(x)\). Compactness of
\(I\times\{s:\delta\le |s|\le\pi\}\) gives a uniform gap
\[
  |t|\ge(1+\eta)\tau(x).
\]
Possible apparent singularities of the closed forms at the zeros of
\(Q\) are removable for \(w(z,t)\), since at a zero \(z_0\) of \(Q\) the
balanced PDE reduces to \(w_t=(\alpha_0z_0+\gamma_0)w\).
\end{proof}

\begin{proof}[Proof of Lemma~\ref{lem:fourier-decay}]
Let \(x=e^\theta\), \(z=xe^{is}\), and \(I=e^\Theta\). By
Lemma~\ref{lem:singularity-separation}, \(w(z,t)\) is analytic for
\[
  |t|\le(1+\eta)\tau(x)
\]
uniformly for \(x\in I\) and \(\delta\le |s|\le\pi\). By compactness,
\(w(z,t)\) is uniformly bounded on these sets. Cauchy's estimate gives
\[
  |P_n(z)|
  \le
  C_1 n!\big((1+\eta)\tau(x)\big)^{-n}
  \le
  C_1'
  n^{1/2}
  \left(\frac{n}{e(1+\eta)\tau(x)}\right)^n.
\]
On the other hand, Theorem~\ref{thm:quadratic-transfer} gives uniformly
for \(x\in I\)
\[
  P_n(x)
  \ge
  C_2
  n^{\nu-\frac12}
  \left(\frac{n}{e\tau(x)}\right)^n
\]
for all sufficiently large \(n\). Hence
\[
  \left|\frac{P_n(z)}{P_n(x)}\right|
  \le
  C_3 n^{1-\nu}(1+\eta)^{-n}.
\]
Absorbing the polynomial factor into the exponential gives the result.
\end{proof}

\subsection{Proof of the uniform Daniels approximation}

\begin{proof}[Proof of Theorem~\ref{thm:daniels-uniform}]
We use the lattice inversion formula \eqref{eq:daniels-contour}. The
proof has three parts: cumulant control, localisation of the real saddle,
and estimation of the Fourier integral.

First, Theorem~\ref{thm:quadratic-transfer}, applied uniformly in a
complex neighbourhood of \(e^{\Theta_\varepsilon}\), gives
\[
  \kappa_n(\theta)
  =
  c_n+n\psi(\theta)+g(\theta)+O(n^{-1}),
\]
where
\[
  \psi(\theta)=-\log\tau(e^\theta),
  \qquad
  g(\theta)=\log\mathcal H(e^\theta),
\]
and \(c_n\) is independent of \(\theta\). By Cauchy's estimates,
\[
  \kappa_n^{(r)}(\theta)
  =
  n\psi^{(r)}(\theta)+g^{(r)}(\theta)+O(n^{-1}),
  \qquad r=1,2,3,4,5,
\]
uniformly for real \(\theta\in\Theta_\varepsilon\).

By Lemma~\ref{lem:tau-integral},
\[
  \psi'(\theta)=F'(\theta)=u(e^\theta),
  \qquad
  \psi''(\theta)=e^\theta u'(e^\theta)>0.
\]
Hence \(\psi''\) is bounded away from zero on
\(\Theta_\varepsilon\), and
\[
  \kappa_n''(\theta)\asymp n
\]
uniformly on \(\Theta_\varepsilon\). Moreover,
\[
  \frac1n\kappa_n'(\theta)=F'(\theta)+O(n^{-1})
\]
uniformly on \(\Theta_\varepsilon\). Since \(F'\) is strictly increasing
and the endpoints of \(\Theta_\varepsilon\) bracket
\([\varepsilon,1-\varepsilon]\), the equation
\[
  \kappa_n'(\theta)=k
\]
has a unique solution
\(\theta_{n,k}\in\Theta_\varepsilon\) for all sufficiently large \(n\)
and all \(k/n\in[\varepsilon,1-\varepsilon]\).

For such \(k\), the coefficient \(w_{n,k}\) is positive for all large
\(n\): take \(k\) up-steps, followed by \(n-k\) level-steps at height
\(k\). The up-step weights are positive because \(\alpha_0>0\), and the
level weight at height \(k\ge1\) is \(Bk+\gamma_0>0\), since \(B>0\).

Use \eqref{eq:daniels-contour} with
\(\theta=\theta_{n,k}\), and put
\[
  V_{n,k}:=\kappa_n''(\theta_{n,k}).
\]
On the central arc \(|s|\le n^{-2/5}\), Taylor expansion gives
\[
  \kappa_n(\theta_{n,k}+is)
  =
  \kappa_n(\theta_{n,k})
  +iks
  -\frac12V_{n,k}s^2
  -\frac{i}{6}\kappa_n^{(3)}(\theta_{n,k})s^3
  +\frac{1}{24}\kappa_n^{(4)}(\theta_{n,k})s^4
  +O(ns^5).
\]
Since \(V_{n,k}\asymp n\), the change of variables
\(y=s\sqrt{V_{n,k}}\) yields
\[
  \frac{1}{2\pi}
  \int_{|s|\le n^{-2/5}}
  \exp\!\left\{
    \kappa_n(\theta_{n,k}+is)-\kappa_n(0)-k(\theta_{n,k}+is)
  \right\}\,ds
\]
\[
  =
  \frac{
    \exp\{\kappa_n(\theta_{n,k})-\kappa_n(0)-k\theta_{n,k}\}
  }{\sqrt{2\pi V_{n,k}}}
  \left(1+O(n^{-1})\right).
\]
Here
\[
  \frac{\kappa_n^{(3)}(\theta_{n,k})}{V_{n,k}^{3/2}}=O(n^{-1/2}),
  \qquad
  \frac{\kappa_n^{(4)}(\theta_{n,k})}{V_{n,k}^{2}}=O(n^{-1}),
\]
the first-order cubic contribution is odd and integrates to zero, and
the square of the cubic contribution is \(O(n^{-1})\).

On the intermediate arc \(n^{-2/5}\le |s|\le\delta\), with \(\delta>0\)
fixed sufficiently small, Taylor expansion and the uniform lower bound
on \(V_{n,k}/n\) give
\[
  \Re\{\kappa_n(\theta_{n,k}+is)-\kappa_n(\theta_{n,k})\}
  \le -c_1ns^2.
\]
This contribution is exponentially small relative to the central
Gaussian contribution.

On the distant arc \(\delta\le |s|\le\pi\),
Lemma~\ref{lem:fourier-decay} gives
\[
  \left|
  \frac{P_n(e^{\theta_{n,k}+is})}
       {P_n(e^{\theta_{n,k}})}
  \right|
  \le C e^{-c_2n}
\]
uniformly in \(k/n\in[\varepsilon,1-\varepsilon]\). Hence the distant arc
is also exponentially negligible. Combining the three estimates gives
\eqref{eq:daniels-uniform}.
\end{proof}

\subsection{The central Gaussian window as a small-tilt limit}

\begin{proof}[Proof of Corollary~\ref{cor:central-llt}]
In the central window, the saddlepoint satisfies
\[
  \theta_{n,k}
  =
  \frac{k-\mu_n}{\kappa_n''(0)}
  +
  O\!\left(\frac{|k-\mu_n|^2}{n^2}\right).
\]
Since
\[
  \kappa_n''(0)=\sigma_n^2\asymp n,
\]
we have \(\theta_{n,k}=O(n^{-1/2})\) when
\(k=\mu_n+O(\sigma_n)\). Expanding the exponent in
\eqref{eq:daniels-uniform} at \(\theta=0\) gives
\[
  \kappa_n(\theta_{n,k})-\kappa_n(0)-k\theta_{n,k}
  =
  -\frac{(k-\mu_n)^2}{2\sigma_n^2}
  +O(n^{-1/2}),
\]
and
\[
  \kappa_n''(\theta_{n,k})
  =
  \sigma_n^2\left(1+O(n^{-1/2})\right).
\]
Substitution into Theorem~\ref{thm:daniels-uniform} gives
\eqref{eq:central-llt}. If \(k-\mu_n=O(1)\), then
\(\theta_{n,k}=O(n^{-1})\), and the same expansion gives an \(O(n^{-1})\)
relative error.
\end{proof}

\subsection{Ray-scale consequences of the saddlepoint formula}

\begin{proof}[Proof of Corollary~\ref{cor:ray-prefactor}]
Let \(u=k/n\in J\). Choose \(\varepsilon\in(0,\tfrac12)\) such that
\(J\subset[\varepsilon,1-\varepsilon]\), and let
\(\Theta_\varepsilon\) be as in Theorem~\ref{thm:daniels-uniform}. By
Lemma~\ref{lem:limit-CGF}, uniformly in a complex neighbourhood of
\(\Theta_\varepsilon\),
\[
  \kappa_n(\theta)-\kappa_n(0)
  =
  nF(\theta)+G(\theta)+O(n^{-1}),
  \qquad
  G(\theta)=\log\frac{\mathcal H(e^\theta)}{\mathcal H(1)}.
\]
Cauchy's estimates also give
\[
  \kappa_n''(\theta)=nF''(\theta)+O(1)
\]
uniformly on \(\Theta_\varepsilon\).

Let \(\theta(u)\) be defined by \(F'(\theta(u))=u\). The finite-\(n\)
saddle satisfies
\[
  \kappa_n'(\theta_{n,k})=k=nu.
\]
Since
\[
  \kappa_n'(\theta)=nF'(\theta)+G'(\theta)+O(n^{-1})
\]
and \(F''\) is bounded away from zero on \(\Theta_\varepsilon\), we have
\[
  \theta_{n,k}=\theta(u)+O(n^{-1}),
\]
uniformly for \(u\in J\).

Now write the exponent in Daniels' formula as
\[
  \Lambda_n
  =
  \kappa_n(\theta_{n,k})-\kappa_n(0)-k\theta_{n,k}.
\]
Using the cumulant expansion,
\[
  \Lambda_n
  =
  n\bigl(F(\theta_{n,k})-u\theta_{n,k}\bigr)
  +G(\theta_{n,k})
  +O(n^{-1}).
\]
The function \(\theta\mapsto F(\theta)-u\theta\) has a critical point at
\(\theta=\theta(u)\). Hence the displacement
\(\theta_{n,k}-\theta(u)=O(n^{-1})\) changes the leading term only by
\(O(n^{-1})\), and therefore
\[
  \Lambda_n
  =
  -nI(u)+G(\theta(u))+O(n^{-1}).
\]
Similarly,
\[
  \kappa_n''(\theta_{n,k})
  =
  nF''(\theta(u))\bigl(1+O(n^{-1})\bigr).
\]
Substitution into Theorem~\ref{thm:daniels-uniform} gives
\[
  p_{n,k}
  =
  \frac{\exp\{G(\theta(u))\}}
       {\sqrt{2\pi nF''(\theta(u))}}
  \exp\{-nI(u)\}
  \bigl(1+O(n^{-1})\bigr),
\]
uniformly for \(u=k/n\in J\), which is \eqref{eq:ray-prefactor}.
\end{proof}

\section{Local asymptotics in the degenerate regimes \(A=0\)}
\label{app:linear-constant}

This appendix records the exact forms and central asymptotic scales in
the two degenerate regimes \(A=0\). These regimes differ from the
quadratic case \(A>0\): no moving algebraic singularity in the \(x\)-plane
is present, and the natural global scale is not the \(n\)-scale used in
the large-deviation theory of Section~\ref{sec:limit}.

Throughout the probabilistic statements below we assume \(\alpha_0>0\),
so that paths can leave height zero. If \(\alpha_0=0\), the terminal
height is identically zero.

\subsection{Constant drift \((A=B=0)\)}

In the constant-drift case the closed form \eqref{eq:A0Beq0} is
exponential-quadratic:
\[
  w(x,t)
  =
  \exp\!\Big((\alpha_0 x+\gamma_0)t+\frac{\alpha_0 C}{2}t^2\Big).
\]
It is convenient to use the Hermite--Kamp\'e de F\'eriet polynomials \cite{AppellKampe1926}
defined by
\[
  \sum_{n\ge0}H_n(X,Y)\frac{t^n}{n!}
  =
  \exp\!\Big(Xt+\frac{Y}{2}t^2\Big).
\]
With
\[
  X(x)=\alpha_0 x+\gamma_0,\qquad Y=\alpha_0 C,
\]
we obtain the exact identities
\begin{equation}
\label{eq:appA-Pn-H}
  P_n(x)=H_n(X(x),Y),
  \qquad
  w_{n,k}
  =
  \binom{n}{k}\alpha_0^k H_{n-k}(\gamma_0,Y).
\end{equation}

\medskip
\noindent\emph{Asymptotics for \(P_n(x)\).}
Assume first that \(Y>0\) and \(X(x)>0\). Cauchy's coefficient formula
gives
\[
  P_n(x)
  =
  \frac{n!}{2\pi i}
  \oint
  \frac{\exp\!\big(X(x)t+\frac{Y}{2}t^2\big)}{t^{n+1}}\,dt.
\]
The positive saddle \(t_*=t_*(x)\) is determined by
\[
  Yt_*^2+X(x)t_*-(n+1)=0.
\]
The positive saddle is exponentially dominant when \(X(x)>0\), and
standard steepest descent yields
\begin{equation}
\label{eq:appA-Pn-asym}
  P_n(x)
  \sim
  \frac{n!}{\sqrt{2\pi\bigl(Y+(n+1)/t_*^2\bigr)}}\,
  \frac{
    \exp\!\bigl(X(x)t_*+\frac{Y}{2}t_*^2\bigr)
  }{t_*^{n+1}}.
\end{equation}

If \(Y=0\), equivalently \(C=0\) or \(\alpha_0=0\), the formula reduces
to the exact expression
\[
  P_n(x)=(\alpha_0 x+\gamma_0)^n.
\]
For \(\alpha_0>0\) and \(\gamma_0>0\), this is the ordinary binomial
case after normalisation.

\medskip
\noindent\emph{Asymptotics for \(w_{n,k}\).}
The exact coefficient formula in \eqref{eq:appA-Pn-H} reduces the local
analysis of \(w_{n,k}\) to that of
\[
  H_m(\gamma_0,Y),
  \qquad m=n-k.
\]
If \(Y>0\) and \(\gamma_0>0\), then
\[
  H_m(\gamma_0,Y)
  =
  \frac{m!}{2\pi i}
  \oint
  \frac{\exp\!\big(\gamma_0 s+\frac{Y}{2}s^2\big)}{s^{m+1}}\,ds,
\]
and the positive saddle \(s_*>0\), determined by
\[
  Ys_*^2+\gamma_0s_*-(m+1)=0,
\]
is exponentially dominant. Hence
\begin{equation}
\label{eq:appA-Hm-asym}
  H_m(\gamma_0,Y)
  \sim
  \frac{m!}{\sqrt{2\pi\bigl(Y+(m+1)/s_*^2\bigr)}}\,
  \frac{
    \exp\!\bigl(\gamma_0s_*+\frac{Y}{2}s_*^2\bigr)
  }{s_*^{m+1}}.
\end{equation}
Combining \eqref{eq:appA-Hm-asym} with \eqref{eq:appA-Pn-H} gives the
corresponding local asymptotics for \(w_{n,k}\).

If \(\gamma_0=0\), the two saddles \(\pm s_*\) have the same modulus and
one has the exact parity constraint
\[
  H_m(0,Y)=0\quad(m\ \text{odd}),
  \qquad
  H_{2j}(0,Y)=\frac{(2j)!}{2^j j!}\,Y^j.
\]
Consequently, in this case
\[
  w_{n,k}=0
  \qquad\text{unless}\qquad
  n-k\equiv0\pmod2.
\]
Equivalently, \(K_n\equiv n\pmod2\). Any local limit statement must then
be interpreted on the admissible span-two sublattice, with the usual
factor \(2\).

\medskip
\noindent\emph{Moments and central window.}
Since \(H_n\) is an Appell family in the variable \(X\),
\[
  \partial_XH_n(X,Y)=nH_{n-1}(X,Y).
\]
At \(x=1\), with
\[
  X_1=\alpha_0+\gamma_0,
\]
this gives the exact moment identities
\begin{equation}
\label{eq:appA-mu-sigma}
  \mu_n
  =
  \alpha_0\,
  \frac{nH_{n-1}(X_1,Y)}{H_n(X_1,Y)},
  \qquad
  \sigma_n^2
  =
  \alpha_0^2
  \frac{n(n-1)H_{n-2}(X_1,Y)}{H_n(X_1,Y)}
  +\mu_n-\mu_n^2.
\end{equation}
When the variance tends to infinity and the span is one, for instance
when \(\gamma_0>0\), a standard saddle-point/admissibility argument
\cite{Bender1973,Hayman1956} gives the Gaussian central window
\[
  p_{n,k}
  =
  \frac{1}{\sqrt{2\pi\sigma_n^2}}
  \exp\!\left(
    -\frac{(k-\mu_n)^2}{2\sigma_n^2}
  \right)(1+o(1))
\]
uniformly for \(k=\mu_n+O(\sigma_n)\). If \(\gamma_0=0\) and \(Y>0\),
the corresponding statement holds on the admissible parity sublattice
with the span-two correction factor.

\subsection{Linear drift \((A=0,\ B>0)\)}

In the linear-drift case the closed form \eqref{eq:A0Bneq0} can be
written as
\[
  w(x,t)
  =
  \exp\!\Big(a t+y(x)(e^{Bt}-1)\Big),
\]
where
\[
  a=\gamma_0-\frac{\alpha_0C}{B},
  \qquad
  y(x)=\frac{\alpha_0}{B}\left(x+\frac{C}{B}\right).
\]
Equivalently,
\[
  w(x,t)=F(t)e^{\Lambda(t)x},
  \qquad
  \Lambda(t)=\frac{\alpha_0}{B}(e^{Bt}-1),
  \qquad
  F(t)=\exp\!\left(\frac{C}{B}\Lambda(t)+a t\right).
\]
We assume \(\alpha_0>0\), so that \(y(x)>0\) for \(x>0\).

\medskip
\noindent\emph{Exact Poisson representation.}
For fixed \(y>0\), define
\[
  F_n(y)
  =
  n![t^n]\exp\!\big(a t+y(e^{Bt}-1)\big).
\]
Then
\[
  F_n(y)
  =
  e^{-y}\sum_{m\ge0}\frac{y^m}{m!}(a+Bm)^n
  =
  \mathbb E\big[(a+BM_y)^n\big],
  \qquad
  M_y\sim\mathrm{Poisson}(y).
\]
Thus
\[
  P_n(x)=F_n(y(x)).
\]
Moreover, if
\[
  y_0=\frac{\alpha_0C}{B^2},
\]
then the coefficients can be recovered exactly as
\[
  w_{n,k}
  =
  \frac{1}{k!}
  \left(\frac{\alpha_0}{B}\right)^k
  F_n^{(k)}(y_0).
\]
This is often a convenient exact representation, although the resulting
closed formulas for \(w_{n,k}\) are not as compact as in the constant
case.

\medskip
\noindent\emph{Asymptotics for \(P_n(x)\).}
Cauchy's formula gives
\[
  P_n(x)
  =
  \frac{n!}{2\pi i}
  \oint
  \frac{\exp\!\big(a t+y(x)(e^{Bt}-1)\big)}{t^{n+1}}\,dt.
\]
The phase is
\[
  \phi(t;x)
  =
  a t+y(x)(e^{Bt}-1)-(n+1)\log t.
\]
It has a unique positive saddle \(t_*=t_*(x)\) satisfying
\[
  \frac{n+1}{t_*}=a+B\,y(x)e^{Bt_*}.
\]
For large \(n\), this saddle satisfies
\[
  t_*(x)
  =
  \frac{1}{B}
  W\!\left(\frac{n}{y(x)}\right)(1+o(1)),
\]
where \(W\) is the Lambert \(W\)-function \cite{Corless1996}. Evaluation
at the saddle gives
\begin{equation}
\label{eq:appA-Pn-lin}
  P_n(x)
  \sim
  \frac{n!}{\sqrt{2\pi\,\phi''(t_*;x)}}\,
  \frac{
    \exp\!\big(a t_*+y(x)(e^{Bt_*}-1)\big)
  }{t_*^{n+1}},
\end{equation}
where
\[
  \phi''(t_*;x)
  =
  B^2y(x)e^{Bt_*}+\frac{n+1}{t_*^2}.
\]

\medskip
\noindent\emph{Moments and central scale.}
Let
\[
  y_s=y(1)=\frac{\alpha_0}{B}\left(1+\frac{C}{B}\right),
  \qquad
  W_s=W\!\left(\frac{n}{y_s}\right),
\]
and set
\[
  \rho=\frac{B}{B+C}.
\]
Differentiating the saddle expansion at \(x=1\) gives
\begin{equation}
\label{eq:appA-mu-sigma-lin}
  \mu_n
  =
  \rho\,\frac{n}{W_s}+O(1),
  \qquad
  \sigma_n^2
  =
  \rho(1-\rho)\frac{n}{W_s}
  +
  \rho^2\frac{n}{W_s^2+W_s}
  +O(1).
\end{equation}
Thus the natural central scale is sublinear. If \(C>0\), then
\(\rho<1\), and
\[
  \mu_n\asymp \frac{n}{\log n},
  \qquad
  \sigma_n^2\asymp \frac{n}{\log n}.
\]
If \(C=0\), then \(\rho=1\), the first variance term vanishes, and
\[
  \mu_n\asymp \frac{n}{\log n},
  \qquad
  \sigma_n^2\asymp \frac{n}{(\log n)^2}.
\]
The central window is Gaussian after the corresponding sublinear
normalisation. In particular, no non-trivial \(n\)-speed
large-deviation principle arises in the linear-drift regime.

\section{Local asymptotics in the quadratic regime ($A>0$)}
\label{app:quadratic}

This appendix spells out Theorem~\ref{thm:quadratic-transfer} in the
three quadratic Pearson regimes. Throughout we assume the balanced case
\[
  \beta_0=C,\qquad A>0,\qquad \alpha_0>0,
\]
and write
\[
  \nu=\frac{\alpha_0}{A}.
\]
The singular times \(\tau(x)\) are those of
Lemma~\ref{lem:tau-integral}. In all three regimes, the general
coefficient estimate is
\begin{equation}
\label{eq:appB-master}
  P_n(x)
  =
  \frac{\sqrt{2\pi}}{\Gamma(\nu)}\,
  \mathcal H(x)\,
  n^{\nu-\frac12}
  \left(\frac{n}{e\,\tau(x)}\right)^n
  \left(1+O(n^{-1})\right),
\end{equation}
locally uniformly for \(x\) in compact subsets of \((0,\infty)\). The
mean and variance follow from
\[
  \chi(x)
  =
  -\frac{\tau'(x)}{\tau(x)}
  =
  \frac{1}{Q(x)\tau(x)},
  \qquad
  u(x)=x\chi(x),
\]
namely
\begin{equation}
\label{eq:appB-moments-general}
  \mu_n=n\chi(1)+O(1),
  \qquad
  \sigma_n^2=n u'(1)+O(1).
\end{equation}
Equivalently,
\[
  \sigma_n^2
  =
  \mu_n
  +
  n\!\left(
    \chi(1)^2-\frac{\tau''(1)}{\tau(1)}
  \right)
  +O(1).
\]
When \(B+C>0\), Lemma~\ref{lem:tau-integral} gives
\(u'(1)>0\), so the leading variance coefficient is positive. In the
boundary case \(B=C=0\), the leading coefficient vanishes and
\(K_n/n\) degenerates at \(1\).

\subsection{Two real roots \((\Delta>0)\)}

Let
\[
  Q(x)=A(x-r_1)(x-r_2),\qquad r_1<r_2,
\]
and put
\[
  d=r_2-r_1>0,\qquad
  L(x)=\log\frac{x-r_1}{x-r_2},\qquad
  S(x)=(x-r_1)(x-r_2).
\]
Then
\[
  \tau(x)=\frac{L(x)}{A d}.
\]
With
\[
  c_0=\alpha_0r_1+\gamma_0,
\]
Theorem~\ref{thm:quadratic-transfer} gives
\[
  \mathcal H(x)=
  \exp\!\big(c_0\tau(x)\big)
  \big(A\tau(x)(x-r_2)\big)^{-\nu}.
\]
Consequently,
\begin{equation}
\label{eq:appB-Pn-DeltaPos}
  P_n(x)
  =
  \frac{\sqrt{2\pi}}{\Gamma(\nu)}\,
  \exp\!\big(c_0\tau(x)\big)
  \big(A\tau(x)(x-r_2)\big)^{-\nu}
  n^{\nu-\frac12}
  \left(\frac{n}{e\,\tau(x)}\right)^n
  \left(1+O(n^{-1})\right).
\end{equation}
Equivalently, since \(A\tau(x)(x-r_2)=L(x)(x-r_2)/d\), the amplitude can
also be written as
\[
  \mathcal H(x)
  =
  \exp\!\big(c_0\tau(x)\big)
  \left(\frac{d}{L(x)(x-r_2)}\right)^\nu.
\]

The logarithmic derivative and the second derivative ratio are
\[
  \chi(x)
  =
  \frac{d}{S(x)L(x)},
\]
and
\[
  \frac{\tau''(x)}{\tau(x)}
  =
  \frac{d\,\big(2x-r_1-r_2\big)}{S(x)^2L(x)}.
\]
Thus, with
\[
  S_1=(1-r_1)(1-r_2),
  \qquad
  L_1=\log\frac{1-r_1}{1-r_2},
\]
one has
\[
  \mu_n
  =
  n\,\frac{d}{S_1L_1}+O(1),
\]
and
\[
  \sigma_n^2
  =
  n\left[
    \frac{d}{S_1L_1}
    +
    \frac{d^2}{S_1^2L_1^2}
    -
    \frac{d(2-r_1-r_2)}{S_1^2L_1}
  \right]
  +O(1).
\]

\subsection{Double root \((\Delta=0)\)}

Let
\[
  Q(x)=A(x-r)^2,
  \qquad
  r=-\frac{B}{2A},
\]
and put
\[
  c_0=\alpha_0r+\gamma_0.
\]
Then
\[
  \tau(x)=\frac{1}{A(x-r)},
  \qquad
  \mathcal H(x)=\exp\!\big(c_0\tau(x)\big).
\]
Thus
\begin{equation}
\label{eq:appB-Pn-DeltaZero}
  P_n(x)
  =
  \frac{\sqrt{2\pi}}{\Gamma(\nu)}\,
  \exp\!\big(c_0\tau(x)\big)
  n^{\nu-\frac12}
  \left(\frac{n}{e\,\tau(x)}\right)^n
  \left(1+O(n^{-1})\right).
\end{equation}
Here
\[
  \chi(x)=\frac{1}{x-r},
  \qquad
  \frac{\tau''(x)}{\tau(x)}=\frac{2}{(x-r)^2}.
\]
Therefore, at \(x=1\), in the non-degenerate case \(r<0\),
\begin{equation}
\label{eq:appB-mu-sigma-DeltaZero}
  \mu_n=\frac{n}{1-r}+O(1),
  \qquad
  \sigma_n^2=\frac{-r}{(1-r)^2}\,n+O(1).
\end{equation}
If \(r=0\), equivalently \(B=C=0\), the leading variance coefficient is
zero and the macroscopic terminal height degenerates at \(K_n/n=1\).

\subsection{Complex conjugate roots \((\Delta<0)\)}

Let
\[
  p=-\frac{B}{2A},
  \qquad
  q=\frac{\sqrt{-\Delta}}{2A}>0,
  \qquad
  Q(x)=A\big((x-p)^2+q^2\big),
\]
and put
\[
  c_0=\alpha_0p+\gamma_0,
  \qquad
  \Theta(x)=\frac{\pi}{2}-\arctan\frac{x-p}{q}.
\]
Then
\[
  \tau(x)=\frac{\Theta(x)}{Aq}.
\]
The amplitude is
\[
  \mathcal H(x)=
  \exp\!\big(c_0\tau(x)\big)
  \big(A\tau(x)\sqrt{(x-p)^2+q^2}\big)^{-\nu}.
\]
Therefore
\begin{equation}
\label{eq:appB-Pn-DeltaNeg}
  P_n(x)
  =
  \frac{\sqrt{2\pi}}{\Gamma(\nu)}\,
  \exp\!\big(c_0\tau(x)\big)
  \big(A\tau(x)\sqrt{(x-p)^2+q^2}\big)^{-\nu}
  n^{\nu-\frac12}
  \left(\frac{n}{e\,\tau(x)}\right)^n
  \left(1+O(n^{-1})\right).
\end{equation}
Equivalently, since \(A\tau(x)=\Theta(x)/q\),
\[
  \mathcal H(x)
  =
  \exp\!\big(c_0\tau(x)\big)
  \left(
    \frac{q}{\Theta(x)\sqrt{(x-p)^2+q^2}}
  \right)^\nu.
\]

Moreover,
\[
  \chi(x)
  =
  \frac{q}{\big((x-p)^2+q^2\big)\Theta(x)},
\]
and
\[
  \frac{\tau''(x)}{\tau(x)}
  =
  \frac{2(x-p)q}
       {\big((x-p)^2+q^2\big)^2\Theta(x)}.
\]
The moment formulas again follow from
\eqref{eq:appB-moments-general}. Since \(\Delta<0\) implies \(C>0\),
this regime is non-degenerate and the leading variance coefficient
\(u'(1)\) is positive.

\bigskip
\noindent\emph{Summary.}
Across the three quadratic regimes the local behaviour is determined by
\(\tau(x)\) and its derivatives. The Pearson geometry---logarithmic when
\(\Delta>0\), rational when \(\Delta=0\), and trigonometric when
\(\Delta<0\)---feeds directly into \(\tau\), then into
\(\chi=-\tau'/\tau\), and hence into the leading mean and variance of
\(K_n\). In the one-dominant-singularity regimes \(B>0\), the finite-\(n\)
point probabilities are governed uniformly in the interior by the tilted
saddlepoint analysis of Section~\ref{sec:daniels}; in the exceptional
complex-root case \(B=0\), the antipodal correction described in
Remark~\ref{rem:B0-interference} may be needed.

\section*{Declaration of competing interest}
The author declares that he has no known competing financial interests
or personal relationships that could have appeared to influence the work
reported in this paper.

\section*{Funding}
This research did not receive any specific grant from funding agencies
in the public, commercial, or not-for-profit sectors.

\end{document}